\documentclass[12pt]{amsart}
\usepackage{amscd,amssymb,epsfig}%,showkeys}
\usepackage[mathscr]{eucal}
\usepackage[english]{babel}

\input xy
\usepackage[all]{xy}

\begin{document}

%\ifx\jpptexloaded\relax\endinput\else\let\jpptexloaded\relax\fi

\newtheoremstyle{mytheorem}% name
  {}%      Space above, empty = `usual value'
  {}%      Space below
  {\slshape}% Body font
  {}%         Indent amount (empty = no indent, \parindent = para indent)
  {\scshape}% Thm head font
  {.}%        Punctuation after head
  { }%     Space after thm head: " " = normal interword space;
        %       \newline = linebreak
  {}% Thm head spec

\newtheoremstyle{mydefinition}% name
  {}%      Space above, empty = `usual value'
  {}%      Space below
  {\upshape}% Body font
  {}%         Indent amount (empty = no indent, \parindent = para indent)
  {\scshape}% Thm head font
  {.}%        Punctuation after thm head
  { }%     Space after thm head: " " = normal interword space;
        %       \newline = linebreak
  {}% Thm head spec

\theoremstyle{mytheorem}
\newtheorem{lemma}{Lemma}[section]
\newtheorem{prop}[lemma]{Proposition}
\newtheorem{prop_intro}{Proposition}
\newtheorem{cor}[lemma]{Corollary}
\newtheorem{cor_intro}[prop_intro]{Corollary}
\newtheorem{thm}[lemma]{Theorem}
\newtheorem{thm_intro}[prop_intro]{Theorem}
\newtheorem*{thm*}{Theorem}
\theoremstyle{mydefinition}
\newtheorem{rem}[lemma]{Remark}
\newtheorem*{claim*}{Claim}
\newtheorem{rem_intro}[prop_intro]{Remark}
\newtheorem{rems_intro}[prop_intro]{Remarks}
\newtheorem*{notation*}{Notation}
\newtheorem*{warning*}{Warning}
\newtheorem{rems}[lemma]{Remarks}
\newtheorem{defi}[lemma]{Definition}
\newtheorem*{defi*}{Definition}
\newtheorem{defi_intro}[prop_intro]{Definition}
\newtheorem{defis}[lemma]{Definitions}
\newtheorem{exo}[lemma]{Example}
\newtheorem{exo_intro}[prop_intro]{Example}
\newtheorem{exos_intro}[prop_intro]{Examples}

\numberwithin{equation}{section}

\newcommand{\bibURL}[1]{{\unskip\nobreak\hfil\penalty50{\tt#1}}}

\def\ti{-\allowhyphens}

\newcommand{\thismonth}{\ifcase\month % case 0 --- impossible!
  \or January\or February\or March\or April\or May\or June%
  \or July\or August\or September\or October\or November%
  \or December\fi}
\newcommand{\thismonthyear}{{\thismonth} {\number\year}}
\newcommand{\thisdaymonthyear}{{\number\day} {\thismonth} {\number\year}}

\newcommand{\Iso}{\operatorname{Iso}}
\newcommand{\esssup}{\operatorname{ess\,sup}}
\newcommand{\supp}{\operatorname{supp}}
\newcommand{\CC}{{\mathbb C}}
\newcommand{\FF}{{\mathbb F}}
\newcommand{\HH}{{\mathbb H}}
\newcommand{\GG}{{\mathbb G}}
\newcommand{\KK}{{\mathbb K}}
\newcommand{\NN}{{\mathbb N}}
\newcommand{\PP}{{\mathbb P}}
\newcommand{\QQ}{{\mathbb Q}}
\newcommand{\RR}{{\mathbb R}}
\renewcommand{\SS}{{\mathbb S}}
\newcommand{\TT}{{\mathbb T}}
\newcommand{\ZZ}{{\mathbb Z}}

\newcommand{\Bb}{{\mathcal B}}
\newcommand{\Cc}{{\mathcal C}}
\newcommand{\Dd}{{\mathcal D}}
\newcommand{\Ff}{{\mathcal F}}
\newcommand{\Gg}{{\mathcal G}}
\newcommand{\Hh}{{\mathcal H}}
\newcommand{\Kk}{{\mathcal K}}
\newcommand{\Ll}{{\mathcal L}}
\newcommand{\Mm}{{\mathcal M}}
\newcommand{\Nn}{{\mathcal N}}
\newcommand{\Oo}{{\mathcal O}}
\newcommand{\Pp}{{\mathcal P}}
\newcommand{\Rr}{{\mathcal R}}
\newcommand{\Ss}{{\mathcal S}}
\newcommand{\Tt}{{\mathcal T}}
\newcommand{\Vv}{{\mathcal V}}
\newcommand{\Xx}{{\mathcal X}}
\newcommand{\Yy}{{\mathcal Y}}
\newcommand{\Zz}{{\mathcal Z}}

\newcommand{\fraka}{{\mathfrak a}}
\newcommand{\frakg}{{\mathfrak g}}
\newcommand{\frakk}{{\mathfrak k}}
\newcommand{\frakp}{{\mathfrak p}}
\newcommand{\fraks}{{\mathfrak s}}
\newcommand{\fraku}{{\mathfrak u}}

\newcommand{\dD}{{\mathbf D}}
\newcommand{\gG}{{\mathbf G}}
\newcommand{\hH}{{\mathbf H}}
\newcommand{\pP}{{\mathbf P}}
\newcommand{\lL}{{\mathbf L}}
\newcommand{\qQ}{{\mathbf Q}}

\newcommand{\<}{\langle}
\renewcommand{\>}{\rangle}

\def\h{{\rm H}}
\def\hb{{\rm H}_{\rm b}}
\def\ehb{{\rm EH}_{\rm b}}
\def\essim{\mathrm{EssIm}}
\def\ha{{\rm H}_{(G,K)}}
\def\hc{{\rm H}_{\rm c}}
\def\hcb{{\rm H}_{\rm cb}}
\def\hdr{{\rm H}_{\rm dR}}
\def\ehbc{{\rm EH}_{\rm cb}}
\def\ind{\mathrm{ind}}
\def\Ind{\mathrm{Ind}}
\def\ltwo{\mathrm{L}^2}
\def\lp{\mathrm{L}^p}
\def\linfty{\mathrm{L}^\infty}
\def\linftyw{\mathrm{L}^\infty_{\rm w*}}
\def\linftya{\mathrm{L}^\infty_{\mathrm{w*,alt}}}
\def\la{\mathrm{L}^\infty_{\mathrm{alt}}}
\def\cb{{\rm C}_{\rm b}}
\def\binfty{\mathcal B^\infty_{\mathrm {alt}}}
\def\one{\mathbf{1\kern-1.6mm 1}}
\def\sous#1#2{{\raisebox{-1.5mm}{$#1$}\backslash \raisebox{.5mm}{$#2$}}}
\def\rest#1{{\raisebox{-.95mm}{$\big|$}\raisebox{-2mm}{$#1$}}}
\def\homeo#1{{\sl H\!omeo}^+\!\left(#1\right)}
\def\thomeo#1{\widetilde{{\sl \!H}\!omeo}^+\!\left(#1\right)}
\def\bu{\bullet}
\def\weak{weak-* }
\def\property{\textbf{\rm\textbf A}}
\def\cont{\mathcal{C}}
\def\id{{\it I\! d}}
\def\opposite{^{\rm op}}
\def\oddex#1#2{\left\{#1\right\}_{o}^{#2}}
\def\comp#1{{\rm C}^{(#1)}}
\def\ro{\varrho}
\def\ti{-\allowhyphens}
\def\lra{\longrightarrow}

\def\adg{\operatorname{ad}_\frakg}
\def\bg{B_\frakg}
\def\creg{C_{\rm reg}}
\def\cs{\check S}
\def\deta{{\operatorname{det}_A}}
\def\fib{(G/Q)^n_f}
\def\gmodp{\gG(\RR)/\pP(\RR)}
\def\gmodq{\gG(\RR)/\qQ(\RR)}
\def\gr{{\rm Gr}_p(V)}
\def\gra{{\rm Gr}^A_p(V)}
\def\Gr{\mathrm{Gr}}
\def\ghtp{generalized Hermitian triple product }
\def\htp{Hermitian triple product }
\def\isp{\operatorname{Is}_{\<\cdot,\cdot\>}}
\def\isptwo{\operatorname{Is}_{\<\cdot,\cdot\>}^{(2)}}
\def\ispth{\operatorname{Is}_{\<\cdot,\cdot\>}^{(3)}}
\def\isf{\operatorname{Is}_F}
\def\isfi{\operatorname{Is}_{F_i}}
\def\isft{\operatorname{Is}_F^{(3)}}
\def\isfit{\operatorname{Is}_{F_i}^{(3)}}
\def\isftwo{\operatorname{Is}_F^{(2)}}
\def\lin{\operatorname{Lin}(L_+,L_-)}
\def\ll{{\Ll_1,\Ll_2}}
\def\kahler{ K\"ahler }
\def\kg{\kappa_G}
\def\kgb{\kappa_G^{\rm b}}
\def\kib{\kappa_i^{\rm b}}
\def\kibt{\tilde\kappa_i^{\rm b}}
\def\oc{\overline c}
\def\om{\overline m}
\def\pii{\pi_i}
\def\piit{\tilde\pi_i}
\def\psupq{{\rm PSU}(p,q)}
\def\psuvi{{\rm PSU}\big(V,\<\cdot,\cdot\>_i\big)}
\def\rg{r_G}
\def\slv{{\rm SL}(V)}
\def\stab{\operatorname{Stab}}
\def\supq{{\rm SU}(p,q)}
\def\SUq{\mathrm{SU}(q,1)}
\def\SUp{\mathrm{SU}(p,1)}
\def\suv{{\rm SU}\big(V,\<\cdot,\cdot\>\big)}
\def\suvi{{\rm SU}\big(V,\<\cdot,\cdot\>_i\big)}
\def\jm{J^{(m)}}
\def\val{\operatorname{Val}}
\def\vc{V_\CC}
\def\vconj{\overline v}
\def\vol{\operatorname{vol}}
\def\wconj{\overline w}
\def\xh{\mathcal X^h}
\def\xb{\mathcal X^b}
\def\xbc{\overline{\mathcal X^b}}
\def\xhc{\overline{\mathcal X^h}}
\def\xbcth{\overline{\mathcal X^b}^{(3)}}

\def\eI{{^{I}\!\mathrm{E}}}
\def\eII{{^{I\!I}\!\mathrm{E}}}
\def\dI{{^{I}\!d}}
\def\dII{{^{I\!I}\!d}}
\def\fI{{^{I}\!F}}
\def\fII{{^{I\!I}\!F}}

\def\No{N\raise4pt\hbox{\tiny o}\kern+.2em}
\def\no{n\raise4pt\hbox{\tiny o}\kern+.2em}
\def\bsl{\backslash}

%% Format de figures
%%%%

\title[Bounded Cohomology and Deformation Rigidity]{Bounded Cohomology and Deformation Rigidity in Complex Hyperbolic Geometry}
\author{Marc Burger}
\email{burger@math.ethz.ch}
\address{FIM, ETH Zentrum, R\"amistrasse 101, CH-8092 Z\"urich, Switzerland}
\author{Alessandra Iozzi}
\email{Alessandra.Iozzi@unibas.ch}
\address{Institut f\"ur Mathematik, Universit\"at Basel, Rheinsprungstrasse 21,
CH-4051 Basel, Switzerland}
\address{D\'epartment de Math\'ematiques, Universit\'e de Strasbourg, 7, rue Ren\'e Descartes, F-67084 Strasbourg Cedex, France}
\thanks{A.I. was partially supported by FNS grant PP002-102765}
%\keywords{ }
%\subjclass{ }

\date{\today}

\begin{abstract}
We develop further basic tools 
in the theory of bounded continuous cohomology 
of locally compact groups;  as such, this paper can be considered 
a sequel to \cite{Burger_Monod_GAFA}, \cite{Monod_book}, and
\cite{Burger_Iozzi_app}.  

We apply these tools to establish a Milnor--Wood type inequality in a very general
context and to prove a global rigidity result 
which was originally announced in \cite{Burger_Iozzi_00} and \cite{Iozzi_ern}
with a sketch of a proof using bounded cohomology techniques 
and then proven by Koziarz and Maubon
in \cite{Koziarz_Maubon} using harmonic map techniques.
As a corollary one obtains that a lattice in $\SUp$ cannot be deformed 
nontrivially in $\SUq$, $q\geq p$, if either $p\geq2$ or the lattice
is cocompact. 
This generalizes to noncocompact lattices a theorem of 
Goldman and Millson, \cite{Goldman_Millson}.
\end{abstract}

\maketitle

\vskip2cm

\setcounter{tocdepth}{3}
\tableofcontents

\vskip2cm
\section{Introduction}\label{sec:intro}
The continuous cohomology $\hc^\bullet(G,\RR)$ of a topological group $G$
is the cohomology of the complex $(\mathrm{C}(G^\bullet)^G,d^\bullet)$ 
of $G$-invariant continuous functions, while its bounded continuous
cohomology $\hcb^\bullet(G,\RR)$ is the cohomology of the subcomplex
$(\mathrm{C}_\mathrm{b}(G^\bullet)^G,d^\bullet)$ of $G$-invariant
bounded continuous functions.  The inclusion of the complex 
of bounded continuous functions into the one consisting of
continuous functions gives rise to the comparison map
\begin{equation*}
{\tt c}_G^\bullet:\hcb^\bullet(G,\RR)\to \hc^\bullet(G,\RR)
\end{equation*}
which encodes subtle properties of $G$ of algebraic and geometric nature,
see \cite{Bavard}, \cite{Ghys_87}, \cite{Mineyev}, \cite{Burger_Monod_JEMS},
\cite[\S~V.13]{Burger_Monod_GAFA}, \cite{Burger_Iozzi_supq}, 
(see also \cite{Brooks}, \cite{Brooks_Series}, \cite{Grigorchuk}, 
\cite{Soma}, \cite{Bestvina_Fujiwara}, \cite{Epstein_Fujiwara}, 
\cite{Gambaudo_Ghys}, \cite{Entov_Polterovich}, \cite{Biran_Entov_Polterovich}, 
\cite{Kotschick_not}, \cite{Kotschick_proc}
in relation with the existence of quasi-morphisms).
We say that a continuous class on $G$ is representable by a bounded 
continuous class if it is in the image of ${\tt c}_G^\bullet$.  

When $G$ is a connected semisimple Lie group with finite center 
and associated symmetric space $\Xx$ and $L<G$ is any closed subgroup, 
a useful tool in the study of the continuous cohomology of $L$ 
is the van~Est isomorphism,
according to which $\hc^\bullet(L,\RR)$ is canonically isomorphic
to the cohomology $\h^\bullet\big(\Omega^\bullet(\Xx)^L\big)$
of the complex $\big(\Omega^\bullet(\Xx)^L,d^\bullet\big)$
of $L$-invariant smooth differential forms $\Omega^\bullet(\Xx)$ on $\Xx$.
For example, if $\Gamma<G$ is a torsionfree discrete subgroup, 
$\h^\bullet(\Gamma,\RR)$ is the de Rham cohomology 
$\hdr^\bullet(\Gamma\backslash \Xx)$ of the manifold $\Gamma\backslash \Xx$.
(Here and in the sequel we drop the subscript $\ _\mathrm{c}$
if the group is discrete.)  For simplicity, in the introduction 
we restrict ourselves to this case, and we refer the reader 
to the body of the paper for the general statement in the case
in which $\Gamma$ is an arbitrary closed subgroup. 

We do not know of an analogue of van Est theorem in the context
of continuous bounded cohomology.  This paper however explores a particular 
aspect of the comparison map and the pullback,
namely the relation between bounded continuous cohomology
and the complex of, loosely speaking, invariant smooth differential
forms with some boundedness condition. For instance, our first result
gives us information on the differential forms that one can use
to represent a class in the image of the comparison map.

\begin{thm_intro}\label{thm:thm1} Let $\Gamma<G$ be a torsionfree discrete subgroup
of a connected semisimple Lie group $G$ with finite center, 
and $\rho:\Gamma\to G'$ a homomorphism into a topological group $G'$.
If $\alpha\in\hc^n(G',\RR)$ is representable by a continuous bounded class,
then its pullback 
$\rho^{(n)}(\alpha)\in\h^n(\Gamma,\RR)\cong\hdr^n(\Gamma\backslash \Xx)$
is representable by a closed differential $n$-form 
on $\Gamma\backslash \Xx$ which is bounded.
\end{thm_intro}

Here a form is bounded on $\Gamma\backslash \Xx$ if its supremum
norm, computed using the Riemannian metric, is finite.
We shall see later that in the case in which $G,G'$ are the connected 
components of the isometry groups of complex hyperbolic spaces
and $\alpha$ is the \kahler class, 
the bounded closed $2$-form in Theorem~\ref{thm:thm1} 
can be given explicitly (see Theorem~\ref{thm:thm9}).  

In particular, by taking in the above theorem $\Gamma=G'$
and $\rho=\id$, we obtain:

\begin{cor_intro}\label{cor:cor2} Let $\Gamma<G$ be a torsionfree discrete subgroup
of a connected semisimple Lie group $G$ with finite center and 
associated symmetric space $\Xx$.  
Any class in the image of the comparison map
\begin{equation*}
{\tt c}_\Gamma^\bullet:
\hb^\bullet(\Gamma,\RR)\to\h^\bullet(\Gamma,\RR)\cong\hdr^\bullet(\Gamma\backslash \Xx)
\end{equation*}
is representable by a closed form on $\Gamma\backslash \Xx$ which is bounded.
\end{cor_intro}

Even if $G'$ is a connected Lie group, little is known about
the surjectivity properties of the comparison map ${\tt c}_{G'}^\bullet$.
However, as a direct consequence of a theorem of Gromov \cite{Gromov_82}
which asserts that characteristic classes are bounded
(see \cite{Bucher_thesis} for a resolution of singularities free proof),
we have the following:

\begin{cor_intro}\label{cor:cor3} Let $\Gamma<G$ be a torsionfree discrete
subgroup of a connected semisimple Lie group with finite center $G$
and $\rho:\Gamma\to G$ a homomorphism into a real algebraic group $G'$.
If $\alpha\in\hc^n(G',\RR)$
comes from a characteristic class of a flat principal $G'$-bundle,
then $\rho^{(n)}(\alpha)\in\h^n(\Gamma,\RR)\cong\hdr^n(\Gamma\backslash \Xx)$
is representable by a closed differential $n$-form on $\Gamma\backslash \Xx$
which is bounded.
\end{cor_intro}

Notice that Theorem~\ref{thm:thm1}, and hence Corollary~\ref{cor:cor2}
and Corollary~\ref{cor:cor3} are valid, with an appropriate
formulation, for any closed subgroup $\Gamma<G$ 
(compare with Corollary~\ref{cor:thm1-L} and 
Proposition~\ref{prop:transfer-L}).

If $L$ is a connected semisimple group with finite center, 
one has full information about the comparison map only 
in degree two, in which case
\begin{equation*}
{\tt c}_L^{(2)}:\hcb^2(L,\RR)\to\hc^2(L,\RR)
\end{equation*}
is an isomorphism, \cite{Burger_Monod_GAFA}.  This is the case we exploit,
also because in this degree continuous cohomology admits a simple description.
Recall in fact that if $\Yy$ is the symmetric space associated to $L$, 
the dimension of $\hc^2(L,\RR)$ is the number of irreducible
factors of $\Yy$ of Hermitian type and $\Omega^2(\Yy)^L$
is generated by the \kahler forms of the irreducible Hermitian factors of $\Yy$.
We say that $L$ is of {\it Hermitian type} if $\Yy$ is Hermitian symmetric
and we denote by $\omega_\Yy$ the \kahler form on $\Yy$ and
by $\kappa_\Yy\in\hc^2(L,\RR)$ the corresponding continuous class
under the isomorphism $\hc^2(L,\RR)\cong\h^2\big(\Omega^\bullet(\Yy)^L\big)$.

Let $G$ be of Hermitian type with associated symmetric space $\Xx$
and $\Gamma<G$ a torsionfree lattice;  for $1\leq p\leq\infty$
let $\h_p^\bullet(\Gamma\backslash \Xx)$ denote the $\lp$-cohomology
of $\Gamma\backslash \Xx$, which is the cohomology of the complex
of smooth differential forms $\alpha$ on $\Gamma\backslash \Xx$
such that $\alpha$ and $d\alpha$ are in $\lp$.  Inclusion
in the complex of smooth differential forms gives thus a comparison map
\begin{equation*}
i_p^\bullet:
 \h_p^\bullet(\Gamma\backslash \Xx)\to\hdr^\bullet(\Gamma\backslash \Xx)\,.
\end{equation*}
Then we have:

\begin{cor_intro}\label{cor:cor4} Assume that $G,G'$ are of Hermitian type,
let $\Gamma<G$ be a torsionfree lattice, $\Xx$ the Hermitian symmetric space
associated to $G$ and $\rho:\Gamma\to G'$
a homomorphism.  Then for every $1\leq p\leq\infty$
there is a linear map
\begin{equation*}
\rho_p^{(2)}:\hc^2(G',\RR)\to\h_p^2(\Gamma\backslash \Xx)
\end{equation*}
such that the diagram
\begin{equation*}
\xymatrix{
 \hc^2(G',\RR)\ar[dr]_{\rho_p^{(2)}}\ar[rr]^{\rho^{(2)}}
&
&\h^2(\Gamma,\RR)\ar[r]^\cong
&\hdr^2(\Gamma\backslash \Xx)\\
&\h_p^2(\Gamma\backslash \Xx)\ar[urr]_{i_p^{(2)}}\\
}
\end{equation*}
commutes.
\end{cor_intro}

In the situation in which $\Gamma<G$ is a lattice and $\Xx$ is Hermitian symmetric,
the $\ltwo$-cohomology $\h_2^\bullet(\Gamma\backslash \Xx)$ is reduced 
({\it i.\,e.} Hausdorff) and finite dimensional in all degrees;
it may hence be identified with the space of $\ltwo$-harmonic forms on
$\Gamma\backslash \Xx$ and carries a natural scalar product 
$\langle\,\cdot\,,\,\cdot\,\rangle$.  The \kahler form $\omega_{\Gamma\backslash \Xx}$
is thus a distinguished element of $\h_2^2(\Gamma\backslash \Xx)$.
Given now a homomorphism $\rho:\Gamma\to G'$ and using Corollary~\ref{cor:cor4},
the invariant 
\begin{equation}\label{eq:ipi-intro} 
i_\rho:=\frac{\langle\rho_2^{(2)}(\kappa_{\Xx'}),\omega_{\Gamma\backslash \Xx}\rangle}
             {\langle\omega_{\Gamma\backslash \Xx},\omega_{\Gamma\backslash \Xx}\rangle}
\end{equation}
is well defined and finite.  We have then the following
Milnor--Wood type inequality:

\begin{thm_intro}\label{prop:prop5} Let $G,G'$ be of Hermitian type
with associated symmetric spaces $\Xx$ and $\Xx'$, let 
$\rho:\Gamma\to G$ be a representation of a lattice in $G$ 
with invariant $i_\rho$ as in (\ref{eq:ipi-intro}).
Assume that $\Xx$ is irreducible
and that the Hermitian metrics on $\Xx$ and $\Xx'$ are normalized
so as to have minimal holomorphic sectional curvature $-1$.
Then
\begin{equation}\label{eq:milnor-wood}
|i_\rho|\leq\frac{\operatorname{rk} \Xx'}{\operatorname{rk} \Xx}\,.
\end{equation}
\end{thm_intro}

Let us call a representation $\rho:\Gamma\to G'$ {\it maximal}
when equality holds in (\ref{eq:milnor-wood}).  We assume this and we distinguish
then the following cases:

\medskip
%\noindent
\underline{$\operatorname{rk} \Xx\geq2$.} Using Margulis' superrigidity
\cite{Margulis_book},
one sees that equality occurs if and only if there exists
an equivariant isometric embedding
\begin{equation*}
f:\Xx\to \Xx'
\end{equation*}
which is furthermore {\it tight}.  
Tightness is an analytic condition
which singles out a certain type of isometric embeddings between
Hermitian symmetric spaces, not necessarily holomorphic, and we refer to
\cite{Wienhard_thesis} and \cite{Burger_Iozzi_Wienhard_tight}
for a systematic study of these. 

\medskip
%\noindent
\underline{$\operatorname{rk} \Xx=1$.} In this case $\Xx=\Hh_\CC^p$ is a
complex hyperbolic space and we remark that one does not know to which
extent lattices in $\SUp$ have superrigidity properties,
at least when $p\geq4$.
We distinguish once again two cases.  

The case in which $p=1$ and $\Xx'$ 
is a general Hermitian symmetric space is the object 
of an ongoing study in \cite{Burger_Iozzi_Wienhard_tol}.
In this situation maximal representations lead to new interesting
Kleinian groups in higher rank.

If of the other hand $p\geq2$, we expect maximal representations to come from tight 
embeddings, but cannot rely on any general superrigidity result.
In this paper, we study the case where $\operatorname{rk} \Xx'=1$ and 
we establish the following: 

\begin{thm_intro}\label{thm:thm6} Let $\Gamma<\SUp$ be a lattice 
and $\rho:\Gamma\to\mathrm{PU}(q,1)$ be a maximal representation.
Assume that $p\geq2$.  Then there is an equivariant
isometric embedding
\begin{equation*}
\varphi:\Hh_\CC^p\to\Hh_\CC^q
\end{equation*}
which is holomorphic if $i_\rho=1$ and antiholomorphic 
if $i_\rho=-1$.
\end{thm_intro}

V.~Koziarz and J.~Maubon gave in \cite{Koziarz_Maubon}
a proof of Theorem~\ref{thm:thm6} using harmonic map techniques.
We refer to the introduction of their article 
for an excellent overview of the history and the context 
of the subject.

Concerning the case $p=1$ we have:

\begin{thm_intro}[\cite{Burger_Iozzi_letter}]\label{thm:thm7b}
Let $\Gamma<\mathrm{SU(1,1)}$ 
be a lattice and $\rho:\Gamma\to\mathrm{PU}(q,1)$ a representation such that
$|i_\rho|=1$.  Then $\rho(\Gamma)$ leaves a complex geodesic
invariant.
\end{thm_intro}

This was proven by Toledo \cite{Toledo_89} if $\Gamma$ is a compact surface
group. In the noncompact case a variant of Theorem~\ref{thm:thm7b} was obtained 
by Koziarz and Maubon in \cite{Koziarz_Maubon}, with another
definition of maximality which probably coincides with ours.
Thus Theorem~\ref{thm:thm7b} reduces the study of maximal
representations into $\mathrm{PU}(q,1)$ to the case $q=1$,
for which we have the following:

\begin{thm_intro}[\cite{Burger_Iozzi_letter}]\label{thm:thm8b}
Let $\Gamma<\mathrm{SU(1,1)}$ 
be a lattice and $\rho:\Gamma\to\mathrm{PU}(1,1)$ a representation such that
$|i_\rho|=1$.  Then $\rho(\Gamma)$ is discrete and, modulo the center of 
$\Gamma$, $\rho$ is injective.  In fact, there is a continuous
surjective map $f:\partial\Hh_\CC^1\to\partial\Hh_\CC^1$ such that:
\begin{enumerate}
\item $f$ is weakly order preserving;
\item $f\big(\rho(\gamma)\xi\big)=\gamma f(\xi)$ for all $\gamma\in\Gamma$
and all $\xi\in\partial\Hh_\CC^n$.
\end{enumerate}
Furthermore, if one of the following two assumptions is 
verified:
\begin{enumerate}
\item[(i)] $\rho(\Gamma)$ is a lattice or  
\item[(ii)] $\rho(\gamma)$ is a parabolic element if $\gamma$ is a parabolic
element,
\end{enumerate}
then $f$ is a homeomorphism and $\rho(\Gamma)$ is a lattice.
\end{thm_intro}

Recall that, in the terminology of \cite{Iozzi_ern} 
a map $f:\partial\Hh_\CC^1\to\partial\Hh_\CC^1$ is {\it weakly order preserving} 
if whenever $\xi,\eta,\zeta\in\partial\Hh_\CC^1$ are distinct points 
such that $f(\xi),f(\eta),f(\zeta)\in\partial\Hh_\CC^1$ are also distinct, 
then the two triples have the same orientation.

\begin{exo_intro}
We give an example that shows that the map $f$ is not necessarily 
a homeomorphism.  To this purpose, let us realize the free group on
two generators in two different ways:
\begin{itemize}
\item Let $\Gamma=\,<a,b>$ be the lattice in $\mathrm{PU}(1,1)$ 
generated by the parabolic elements $a$ and $b$ with fundamental
domain an ideal triangle.
\begin{figure}
\includegraphics[scale=0.497]{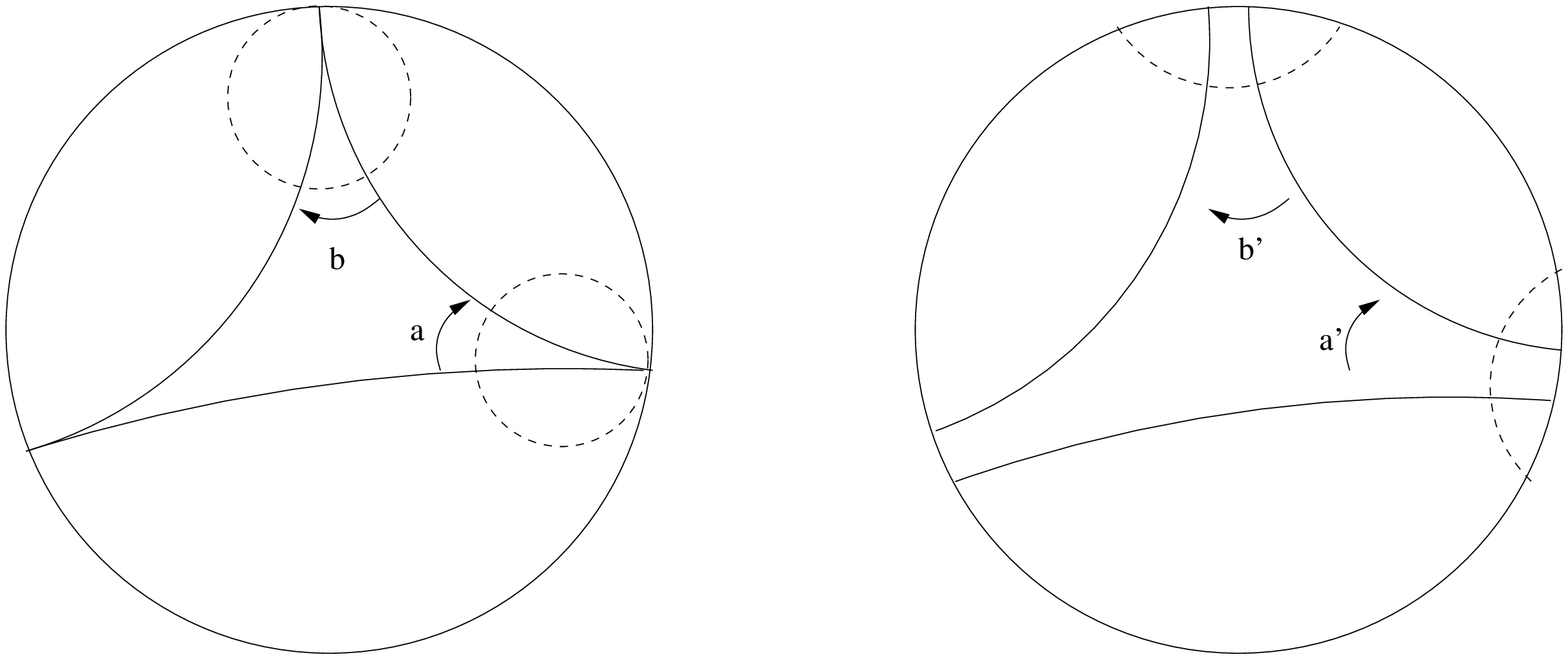}
 \caption{$\Gamma=\,<a,b>$ and $\Lambda=\,<a',b'>$}
\end{figure}
\item Let $\Lambda=\,<a',b'>$ be the convex cocompact group
generated by the hyperbolic elements $a'$ and $b'$.
\end{itemize}

Let $\rho:\Gamma\to\Lambda$ be the representation defined by $\rho(a)=a'$
and $\rho(b)=b'$.  Since $\Lambda$ acts convex cocompactly on $\Hh_\CC^1$,
the orbit map $\Lambda\to\Lambda x$, for $x\in\Hh_\CC^1$ is a quasi-isometry
which extends to a homeomorphism $f_\Lambda:\partial\FF_2\to\Ll_\Lambda$,
where $\FF_2$ is the free group in two generators and $\Ll_\Lambda$ is
the limit set of $\Lambda$.  Likewise, the orbit map
$\Gamma\to\Gamma x$ extends to a continuous surjective map
$f_\Gamma:\partial\FF_2\to\SS^1$ which is one-to-one except for 
the cusps of $\Gamma$, where it is two-to-one.
Then $f_\Gamma\circ f_\Lambda^{-1}$ extends to a map 
$f:\partial\Hh_\CC^1\to\partial\Hh_\CC^1$ such that 
\begin{enumerate}
\item $f$ is weakly order preserving, and
\item $f\big(\rho(\gamma)\xi\big)=\gamma f(\xi)$,
\end{enumerate}

One can prove, using the results in \cite{Ghys_87}, \S~\ref{sec:Toledo}
and \S~\ref{subsec:hermitian} that $i_\rho=1$.
\end{exo_intro}

\medskip

We turn now to the description of the tools involved in the proof
of Theorem~\ref{thm:thm6}.  The ideal boundary $\partial\Hh_\CC^\ell$
of complex hyperbolic $\ell$-space $\Hh_\CC^\ell$ carries a rich geometry 
whose ``lines''are the {\it chains}.  A chain in $\partial\Hh_\CC^\ell$
is by definition the boundary of a complex geodesic in $\Hh_\CC^\ell$;
as such it is a circle equipped with a canonical orientation,
and it is uniquely determined by any two points lying on it.
The ``geometry of chains'' was first studied by E.~Cartan 
who showed that, analogously to the Fundamental Theorem of Projective
Geometry \cite[Theorem~2.26]{Artin}, 
any automorphism of the incidence graph of the geometry of chains comes,
for $\ell\geq2$, from an isometry of $\Hh_\CC^\ell$, \cite{Cartan}.
Closely connected to this is Cartan's {\it invariant angulaire}
\begin{equation*}
c_\ell:(\partial\Hh_\CC^\ell)^3\to[-1,1]
\end{equation*}
introduced in the same paper \cite{Cartan}, which is a full invariant
for the $\mathrm{SU}(\ell,1)$-action on triples of points
in $\partial\Hh_\CC^\ell$ assuming its maximum modulus
(namely 1) exactly on triples of points lying on a chain.

Now let $\Gamma<\SUp$ be a lattice, $\rho:\Gamma\to\mathrm{PU}(q,1)$ 
a homomorphism with nonelementary image
and let $\varphi:\partial\Hh_\CC^p\to\partial\Hh_\CC^q$ be the
$\Gamma$-equivariant measurable boundary map 
whose existence is recalled in \S~\ref{sec:formula}.
A sizable portion of this paper is devoted to developing certain
general tools in the theory of continuous bounded cohomology 
with appropriate coefficients, which will serve to establish a concrete link
between our invariant $i_\rho$ associated to a homomorphism
$\rho:\Gamma\to\mathrm{PU}(q,1)$, the properties
of the corresponding $\Gamma$-equivariant measurable boundary map
$\varphi:\partial\Hh_\CC^p\to\partial\Hh_\CC^q$, and Cartan's
angular invariant.  To illustrate this, let $\Cc_p$ be the set
of chains in $\partial\Hh_\CC^p$, and for almost every chain
$C\in\Cc_p$ (in a sense made precise in \S~\ref{sec:cartan}), 
let us denote by $\varphi_C$ the restriction of $\varphi$ to $C$.
Denoting by $\mu$ the $\SUp$-invariant probability measure on 
$\Gamma\backslash\SUp$, we establish the following formula
which gives a measure of how much the boundary map $\varphi$
distorts a typical chain:

\begin{thm_intro}\label{thm:thm8} Let $\Gamma<\SUp$ be a lattice,
$\rho:\Gamma\to\mathrm{PU}(q,1)$ a homomorphism with nonelementary image and 
$\varphi:\partial\Hh_\CC^p\to\partial\Hh_\CC^q$ be the associated
$\Gamma$-equivariant measurable boundary map.
For almost every chain $C\in\Cc_p$ and
almost every triple $(\xi,\eta,\zeta)\in C^3$, we have
\begin{equation*}
 \int_{\Gamma\backslash\SUp} 
   c_q\big(\varphi_C(g\xi),\varphi_C(g\eta),\varphi_C(g\zeta)\big)\,d\mu(\dot g)
=i_\rho c_p(\xi,\eta,\zeta)\,,
\end{equation*}
where $i_\rho$ is defined in (\ref{eq:ipi-intro}), 
and $c_q$, $c_p$ are the Cartan invariants.
\end{thm_intro}

In the same vein, if $\xi\in\partial\Hh_\CC^p$ and $x\in\Hh_\CC^p$, 
let $e^\xi(x):=e^{h\beta_\xi(0,x)}$,
where $h$ is the volume entropy of $\Hh_\CC^p$,
$\beta_\xi(0,x)$ is the Busemann function relative to a basepoint $0\in\Hh_\CC^p$, 
and $\mu_0$ the $K=\operatorname{Stab}_{\SUp}(0)$-invariant
probability measure on $\partial\Hh_\CC^p$.  The following
is the more precise form of Theorem~\ref{thm:thm1}
announced above:

\begin{thm_intro}\label{thm:thm9} Let $\Gamma<\SUp$ be any torsionfree
discrete subgroup, $\rho:\Gamma\to\mathrm{PU}(q,1)$ a nonelementary homomorphism
and let $\kappa_q$ denote the \kahler class of $\mathrm{PU}(q,1)$.
The $2$-form 
\begin{equation*}
\int_{(\partial\Hh_\CC^p)^3} 
  e^{\xi_0}\wedge de^{\xi_1}\wedge de^{\xi_2} 
  c_q\big(\varphi(\xi_0),\varphi(\xi_1),\varphi(\xi_2)\big)
  d\mu_0(\xi_0)d\mu_0(\xi_1)d\mu_0(\xi_2)
\end{equation*}
is $\Gamma$-invariant, closed, bounded and represents 
$\rho^{(2)}(\kappa_q)\in\hdr^2(\Gamma\backslash\Hh_\CC^p)$.
\end{thm_intro}

Our rigidity theorem on maximal representations follows from
the formula in Theorem~\ref{thm:thm8} and the following
measurable analogue of Cartan's theorem, to which we alluded above:

\begin{thm_intro}\label{thm:thm10} Let $n\geq2$ and let 
$\varphi:\partial\Hh_\CC^p\to\partial\Hh_\CC^q$
be a measurable map such that:
\begin{enumerate}
\item[(i)] for almost every chain $C$ and almost every triple $(\xi,\eta,\zeta)$
of distinct points on $C$, the triple 
$\varphi(\xi), \varphi(\eta), \varphi(\zeta)$ consists also of distinct
points which lie on a chain and have the same orientation as $(\xi,\eta,\zeta)$;
\item[(ii)] for almost every triple of points $\xi,\eta,\zeta$ not on a chain, 
$\varphi(\xi)$, $\varphi(\eta)$, $\varphi(\zeta)$ are also not on a chain.
\end{enumerate}
Then there is an isometric holomorphic embedding $F:\Hh_\CC^p\to\Hh_\CC^q$
such that $\partial F$ coincides with $\varphi$ almost everywhere.
\end{thm_intro}

\medskip

A different way of looking at $i_\rho$ as a {\it foliated Toledo
number} was suggested to us by F.~Labourie, and goes as follows.
The space of configurations of points lying on chains 
can be seen as the space at infinity of the space of configurations 
of points lying on complex geodesics $\Gg_p$
\begin{equation*}
\Gg\Hh_\CC^p=\bigg\{(x,Y):Y\hbox{ is a complex geodesic and }
    \,x\in Y\subset\Hh_\CC^p\bigg\}
\end{equation*}
which is the total space of a foliation whose leaves are the fibers of the map
\begin{equation*}
\begin{aligned}
\pi_2:\Gg\Hh_\CC^p&\to\Gg_p\\
(x,Y)&\mapsto Y
\end{aligned}
\end{equation*}
which are transverse to the fibers of 
\begin{equation*}
\pi_1:\Gg\Hh_\CC^p\to\Hh_\CC^p
\end{equation*}
which, incidentally, are compact.
Given now $\Gamma<\mathrm{PU}(p,1)$ a torsionfree lattice, 
since $\pi_1$ is $\Gamma$-equivariant,
we get a foliated space $\pi_1:\Gg M\to M$ lying above 
$M=\Gamma\backslash \Hh_\CC^p$, 
where $\Gg M=\Gamma\backslash\Gg\Hh_\CC^p$ is foliated by complex geodesics.
The restriction to the complex geodesics of the pullback
$\pi_1^\ast(\omega_M)$ of the \kahler form $\omega_M$ of $M$,
defines a tangential form $\omega_{\Gg M}$.
If then $\rho:\Gamma\to\mathrm{PU}(q,1)$ is a homomorphism
and $\omega'_\rho$ is a bounded closed representative of the class
$\rho^{(2)}(\kappa_q)\in\hdr^2(M)$, then the tangential form $\Omega'_\rho$,
obtained by restricting $\pi^\ast(\omega'_\rho)$ to the leaves
of the foliations, differs from
$\omega_{\Gg M}$ by a bounded function, whose integral
over $\Gg M$ gives $i_\rho$.

\bigskip

\noindent
{\it Application to Deformation Rigidity. } 
If $\Gamma$ is a discrete finitely generated group and $L$ is a topological group,
the space of homomorphisms $\operatorname{Rep}(\Gamma,L)$
of $\Gamma$ into $L$ is topologized naturally as a closed subset of
$L^S$, where $S$ is a finite generating set of $\Gamma$.
Let $BL$ be the classifying space of continuous principal $L$-bundles,
and $c\in\h^\bullet(BL,\RR)$ a characteristic class.  
It is a standard observation that the map
\begin{equation*}
\begin{aligned}
\operatorname{Rep}(\Gamma,L)&\to\h^\bullet(\Gamma,\RR)\\
\rho\,\,\,\,\,\quad&\mapsto\,\,\,\rho_B^\bullet(c)\,,
\end{aligned}
\end{equation*}
where $\rho^\bullet_B:\h^\bullet(BL,\RR)\to\h^\bullet(B\Gamma,\RR)=\h^\bullet(\Gamma,\RR)$
denotes the pullback, is constant on connected components
of $\operatorname{Rep}(\Gamma,L)$.

Assume now that $L$ is irreducible of Hermitian type and
let $K$ be a maximal compact subgroup of $L$.  
It follows from the Iwasawa decomposition that $BK$ is homotopic
equivalent to $BL$, and by Chern-Weil theory
$\h^\bullet(K,\RR)$ is described by the $K$-invariant
polynomials on the Lie algebra $\mathfrak k$ of $K$.
Since $L$ is irreducible Hermitian, the center $Z(\mathfrak k)$ 
is one dimensional and the orthogonal projection of $\mathfrak k$
on $Z(\mathfrak k)$ gives rise to an invariant linear form which,
via Chern-Weil theory, gives rise to a class in 
$\h^2(BK,\RR)=\h^2(BL,\RR)$.
This class corresponds then via the natural homomorphism 
$\h^2(BL,\RR)\to\hc^2(L,\RR)$ to the \kahler class $\kappa_\Yy$,
and hence the commutativity of the diagram
\begin{equation*}
\xymatrix{
 \h^2(BL,\RR)\ar[r]^{\rho^{(2)}_B}\ar[d]
&\h^2(B\Gamma,\RR)\ar[d]^=\\
 \hc^2(L,\RR)\ar[r]^{\rho^{(2)}}
&\h^2(\Gamma,\RR)
}
\end{equation*}
implies that the map
\begin{equation*}
\begin{aligned}
\operatorname{Rep}(\Gamma,L)&\to\h^2(\Gamma,\RR)\\
\rho\qquad&\mapsto\rho^{(2)}(\kappa_\Yy)\,,
\end{aligned}
\end{equation*}
where $\Yy$ is the symmetric space associated to $L$, 
is constant on connected components of $\operatorname{Rep}(\Gamma,L)$.

To turn to our immediate application, let us assume that $p\leq q$ 
and let $\rho_0:\SUp\to\mathrm{PU}(q,1)$ be a {\it standard} homomorphism,
that is a homomorphism associated to any isometric holomorphic embedding 
\begin{equation*}
F:\Hh_\CC^p\to\Hh_\CC^q\,.
\end{equation*}
Observe that %not only $\rho_0$ is maximal\marginpar{what does maximal mean here?} 
%since the embedding is holomorphic,
any two such embeddings $\Hh_\CC^p\to\Hh_\CC^q$ are conjugate in
$\mathrm{PU}(q,1)$;  moreover, the stabilizer in $\mathrm{PU}(q,1)$ of the image 
of $F$ is the almost direct product of the image $\rho_0(\SUp)$
and its centralizer $Z(\rho_0)$ in $\mathrm{PU}(q,1)$,
which is compact.

\begin{cor_intro}\label{cor:cor10} Let $\rho_0:\SUp\to\mathrm{PU}(q,1)$
be a standard representation, 
let $\Gamma<\SUp$ be a lattice and assume that $p\geq2$.
Then any representation $\rho:\Gamma\to\mathrm{PU}(q,1)$ in the
path connected component of $\rho_0|_\Gamma$ in the representation
variety $\operatorname{Rep}(\Gamma,\mathrm{PU}(q,1))$ is, 
modulo conjugation by $\mathrm{PU}(q,1)$, 
of the form $\rho_0\times\omega$,
where $\omega$ is a homomorphism of $\Gamma$ into the compact 
group $Z(\rho_0)$.  
\end{cor_intro}

\begin{rem_intro} We recall that if $\Gamma<\SUp$ is cocompact,
this was proven by Goldman and Millson in \cite{Goldman_Millson}.  
On the other hand, Gusevskii and Parker found quasi-Fuchsian deformations
of a noncocompact lattice $\Gamma<\mathrm{SU}(1,1)$ into $\mathrm{PU}(2,1)$,
\cite{Gusevskii_Parker}.  %Note that this case falls nicely out of the realm 
%of consideration of our methods since $\h^2(\Gamma,\RR)=0$ for any 
%(torsionfree) noncocompact lattice $\Gamma<\mathrm{SU}(1,1)$.
\end{rem_intro}

\bigskip

\noindent
{\it Organization of the Paper: }
Theorem~\ref{thm:thm1} is proven as Corollary~\ref{cor:thm1-L}, 
Corollary~\ref{cor:cor2} is proven as Proposition~\ref{prop:transfer-L},
Corollary~\ref{cor:cor4} is proven as Corollary~\ref{cor:factor-lp}, 
Theorem~\ref{prop:prop5} follows from Lemma~\ref{lem:milnor-wood}
and Lemma~\ref{lem:eq}, Theorems~\ref{thm:thm6}, ~\ref{thm:thm7b} and 
\ref{thm:thm8b} are proven in \S~\ref{sec:proof},
Theorem~\ref{thm:thm8} is Theorem~\ref{thm:formula}, Theorem~\ref{thm:thm9}
is Proposition~\ref{prop:bounded-rep} and finally Theorem~\ref{thm:thm10} is proven
as Theorem~\ref{thm:cartan}.

\bigskip
\noindent
{\it Acknowledgments:} The authors thank Leslie Saper for 
helpful comments concerning $\L^2$-cohomology,
Theo B\"uhler for detailed comments,
and Fran\c cois Labourie for enlightening conversations. 
%%% Local Variables: 
%%% mode: latex
%%% TeX-master: "deform"
%%% End: 

\vskip1cm
\section{Preliminaries on Bounded Cohomology and Hermitian Symmetric Spaces}
\label{sec:prelim}

Here we collect a few basic facts and definitions 
which will be useful in the sequel.  
Basic references are \cite{Monod_book} and \cite{Goldman_book}.

\bigskip

Let $G$ be a locally compact group.  A {\it coefficient $G$-module} is 
a Banach space $E$ with an isometric $G$-action $\pi:G\to\Iso(E)$ 
contragredient to some separable continuous Banach $G$-module;
in particular the $G$-action on $E$ is continuous in the
w$^\ast$-topology.   Given a coefficient $G$-module $E$,
the {\it bounded continuous cohomology with coefficients in $E$} 
is the cohomology of the complex $(\cb(G^\bullet,E)^G,d^\bullet)$
of the space of continuous bounded maps $G^{n+1}\to E$
which are $G$-invariant with respect to the $G$-action on 
$\cb(G^{n+1},E)$ defined by 
\begin{equation*}
(gf)(x):=\pi(g)^{-1}f(gx)\,,
\end{equation*}
for all $x\in G^{n+1}$ and $g\in G$.
Basic examples of coefficient modules are continuous unitary
representations on a separable Hilbert space $\Hh$-- among which 
the trivial one with $\Hh=\RR$ -- and, more importantly in this paper, 
the space $\linfty(G/H)$, where $H\leq G$ is a closed subgroup.
Notice that $\hcb^n(G,E)$ comes naturally equipped with a seminorm
induced by the supremum norm on $\cb(G^\bullet,E)$ and in 
some cases, as for instance if $n=2$ and the coefficient module is separable,
the seminorm is actually a norm.

There is a notion of {\it relatively injective} $G$-module
which serves for the homological algebra characterization of cohomology.
For the precise definition see \cite{Burger_Monod_GAFA} or \cite{Monod_book},
while for our purpose it will suffice to say that if $(S,\nu)$
is a regular measure $G$-space, then the $G$-module
$\linfty(S)$ is relatively injective if and only if 
the $G$-action on $S$ is amenable in the sense of Zimmer, \cite{Zimmer_book}.

Let $E$ be a coefficient $G$-module and $(E_\bullet,d^\bullet)$ be a complex,
where $E_0=E$ is a coefficient $G$-module and $E_n$, for $n\geq1$
are Banach $G$-modules.  We say that $(E_\bullet,d^\bullet)$
is a {\it strong resolution of $E$} if there is a sequence
$h_n:E_n\to E_{n-1}$ of homotopy operators, such that
\begin{enumerate}
\item[(i)] $\|h_n\|\leq1$;
\item[(ii)] $h_n$ maps the subspace of $G$-continuous vectors $\Cc E_n$
into $\Cc E_{n-1}$.  
\end{enumerate}
Then, if the $E_n$, $n\geq1$, are relatively
injective, the cohomology of the subcomplex
\begin{equation*}
\xymatrix@1{
 0\ar[r]
&E^G\ar[r]
&E_1^G\ar[r]
&\cdots\ar[r]
&E_n^G\ar[r]^{d_n}
&\cdots
}
\end{equation*}
is canonically isomorphic to the
bounded continuous cohomology $\hcb^\bullet(G,E)$.
The following lemma collects many functoriality statements needed
in this paper, and is a small modification of a lemma in \cite{Monod_book}.

\begin{lemma}\label{lem:coeff-gp} Let $G,G'$ be locally compact groups,
$\rho:G\to G'$ a continuous homomorphism,
$E$ a $G$-coefficient module and $F$ a $G'$-coefficient module.
Let $\alpha:F\to E$ be a morphism of $G$-coefficient modules,
where the $G$-module structure on $F$ is via $\rho$.
Let $(E_\bullet)$ be a strong $G$-resolution of $E$ by 
relatively injective $G$-modules, and let $(F_\bullet)$ 
be a strong $G'$-resolution of $F$.  Then any two
extensions of the morphism $\alpha$ to a morphism of
$G$-complexes induce the same map in cohomology
\begin{equation*}
\h^\bullet\big(F_\bullet^{G'}\big)\to\h^\bullet\big(E_\bullet^{G}\big)\,.
\end{equation*}
\end{lemma}

\begin{proof}
By \cite[Lemma~7.2.6]{Monod_book} any two extensions of $\alpha$ are
$G$-homotopic and hence induce the same map in cohomology
\begin{equation*}
\h^\bullet\big(F_\bullet^{\rho(G)}\big)\to\h^\bullet\big(E_\bullet^{G}\big)\,.
\end{equation*}
Moreover, the inclusion of complexes $F_\bullet^{G'}\subset F_\bullet^{\rho(G)}$
induces a unique map in cohomology
\begin{equation*}
\h^\bullet\big(F_\bullet^{G'}\big)\to\h^\bullet\big(F_\bullet^{\rho(G)}\big)\,,
\end{equation*}
hence proving the lemma.
\end{proof}

\bigskip

Let now $G$ be a connected, semisimple Lie group with finite center,
and $\Xx$ the associated symmetric space.  Assume that $\Xx$
is Hermitian symmetric, so that on $\Xx$ there exists
a nonzero $G$-invariant (closed) differential $2$-form, 
namely the \kahler form of the Hermitian metric, 
which we denote by $\omega_\Xx\in\Omega^2(\Xx)^G$.
If $x\in\Xx$ is a reference point, and
$\Delta(g_1x,g_2x,g_3x)\subset\Xx$ is a triangle with vertices 
$g_1x,g_2x,g_3x$, geodesic sides and arbitrarily $\mathrm{C}^1$-filled, 
\begin{equation*}
c(g_1,g_2,g_3):=\int_{\Delta(g_1x,g_2x,g_3x)}\omega_\Xx
\end{equation*}
is a differentiable homogeneous $G$-invariant cocycle
and defines the continuous class $\kappa_\Xx\in\hc^2(G,\RR)$
corresponding to $\omega_\Xx$ by the van Est isomorphism
$\hc^2(G,\RR)\simeq\Omega^2(\Xx)^G$.  Moreover, 
$c$ is bounded (\cite{Domic_Toledo}, \cite{Clerc_Orsted_2}),
and hence it defines a bounded continuous class 
$\kappa_\Xx^\mathrm{b}\in\hcb^2(G,\RR)$ which corresponds
to $\kappa_\Xx\in\hc^2(G,\RR)$ under the isomorphism
\begin{equation*}
\hcb^2(G,\RR)\cong\hc^2(G,\RR)\,.
\end{equation*}
It follows again from \cite{Domic_Toledo} and \cite{Clerc_Orsted_2}
that if we normalize the metric on $\Xx$ so as its minimal holomorphic
sectional curvature is $-1$, 
the Gromov norm of the class $\kappa_\Xx^\mathrm{b}$ is
\begin{equation}\label{eq:gromov-norm}
\|\kappa_\Xx^\mathrm{b}\|=\pi\operatorname{rk}\Xx\,.
\end{equation}
In the special case of complex hyperbolic space $\Hh_\CC^\ell$,
the multiple $\frac{1}{\pi}\kappa_\ell^\mathrm{b}$ of the
bounded \kahler class $\kappa_\ell^\mathrm{b}$ (which is here 
and in the following a shortcut for $\kappa_{\Hh_\CC^\ell}^\mathrm{b}$)
admits an explicit representative on $\partial\Hh_\CC^\ell$ given 
by the Cartan cocycle 
\begin{equation*}
c_\ell:(\partial\Hh_\CC^\ell)^3\to[-1,1],
\end{equation*}
which is defined in terms of the Hermitian triple product of a triple
of points in the underlying complex vector space $V$ of dimension $\ell+1$
with a Hermitian form of signature $(\ell,1)$ whose cone of negative
lines gives a model of complex hyperbolic space $\Hh_\CC^\ell$.
Any $(k+1)$-dimensional nondegenerate indefinite linear subspace 
$W\subset V$ gives rise to a {\it $k$-plane}, 
that is a totally geodesic holomorphically embedded 
isometric copy of $\Hh_\CC^k$, whose boundary 
in $\partial\Hh_\CC^\ell$ is called a {\it $k$-chain}.
In particular, $1$-chains (or {\it chains}) 
are the boundary of complex geodesics and play a fundamental role here.  
We refer the reader to \cite{Goldman_book}
for the precise definitions, and we limit ourselves to recall the
following essential lemma:
\begin{lemma} The Cartan cocycle
$c_\ell:(\partial\Hh_\CC^\ell)^3\to[-1,1]$ is a strict 
$\mathrm{SU(\ell,1)}$-invariant Borel cocycle
and $|c_\ell(a,b,c)|=1$ if and only if $a,b,c$ 
are on a chain and pairwise distinct.
\end{lemma}

%%% Local Variables: 
%%% mode: latex
%%% TeX-master: "deform"
%%% End: 

\vskip1cm
\section{Some Cohomological Tools for Locally Compact Groups}
In this section we develop some tools in bounded cohomology
for locally compact groups and their closed subgroups
which will be applied to our specific situation.  
More precisely, while the functorial machinery developed 
in \cite{Burger_Monod_GAFA}, \cite{Monod_book} and \cite{Burger_Iozzi_app}
applies in theory to general strong resolutions,
in practice one ends up working mostly with spaces
of functions on Cartesian products.  In this section 
we deal with spaces of functions on fibered products 
(of homogeneous spaces), whose general framework would be 
that of complexes of functions on appropriate sequences 
$(\Ss_n,\nu_n)$ of (amenable) spaces which are 
analogues of simplicial sets in the category of measured spaces. 

\subsection{Cohomology with Coefficients: With the Use of Fibered Products}
The invariants we consider in this paper are bounded classes
with trivial coefficients; however applying a judicious 
change of coefficients -- from $\RR$ to the $\linfty$ functions
on a homogeneous space -- 
we capture information which otherwise would be lost 
by the use of measurable maps (see the less cryptic Remark~\ref{rem:fibered}).

\subsubsection{Realization on Fibered Products}\label{subsec:fibered} 
The goal of this section is to define the fibered product 
of homogeneous spaces and 
prove that the complex of $\linfty$ functions on fibered products
satisfy all properties necessary to be used to compute
bounded cohomology.  Observe that because of the projection in 
(\ref{eq:p_n}), we shall deal here with cohomology with coefficients.

\bigskip

Let $G$ be a locally compact, second countable group 
and $Q,H$ closed subgroups of $G$ such that $Q\leq H$.
We define the $n$-fold fibered product $\fib$ of $G/Q$
with respect to the canonical projection $p:G/Q\to G/H$ 
to be, for $n\geq1$, the closed subset of $(G/Q)^n$ defined by 
\begin{equation*}
\fib:=\big\{(x_1,\dots,x_n)\in(G/Q)^n\,:\,\,p(x_1)=\cdots=p(x_n)\big\}\,,
\end{equation*}
and we set $\fib=G/H$ if $n=0$.
The invariance of $\fib$ for the diagonal $G$-action on $(G/Q)^n$
induces a $G$-equivariant projection 
\begin{equation}\label{eq:p_n}
p_n:\fib\to G/H
\end{equation}
whose typical fiber is homeomorphic to $(H/Q)^n$.

A useful description of $\fib$ as a quotient space may be obtained as follows.
Considering $H/Q$ as a subset of $G/Q$, the map
\begin{equation}\label{eq:homeo-fib}
\begin{aligned}
q_n:G\times(H/Q)^n&\to\quad\fib\\
(g,x_1,\dots,x_n)&\mapsto(gx_1,\dots,gx_n)
\end{aligned}
\end{equation}
is well defined, surjective, $G$-equivariant 
(with respect to the $G$-action on the first coordinate on $G\times(H/Q)^n$ 
and the product action on $\fib$) 
and invariant under the right
$H$-action on $G\times(H/Q)^n$ defined by 
\begin{equation}\label{eq:Haction}
(g,x_1,\dots,x_n)h:=(gh,h^{-1}x_1,\dots,h^{-1}x_n)\,.
\end{equation}
It is then easy to see that $q_n$ induces a $G$-equivariant homeomorphism
\begin{equation*}
\big(G\times(H/Q)^n\big)/H\to\fib\,,
\end{equation*}
which hence realizes the fibered product $\fib$ as a quotient space.

Let now $\mu$ and $\nu$ be Borel probability measures respectively on $G$ and $H/Q$,
such that $\mu$ is in the class of the Haar measure on $G$ and 
$\nu$ is in the $H$-invariant measure class on $H/Q$.
The pushforward $\nu_n=(q_n)_*(\mu\times\nu^n)$ of the probability measure 
$\mu\times\nu^n$ under $q_n$
is then a Borel probability measure on $\fib$ whose class is $G$-invariant
and thus gives rise to Banach $G$-modules $\linfty(\fib)$
and $G$-equivariant (norm) continuous maps
\begin{equation*}
d_n:\linfty\big(\fib\big)\to\linfty\big((G/Q)^{n+1}_f\big)\,,\hfill\hbox{for }n\geq0\,,
\end{equation*}
defined as follows:
\begin{enumerate}
\item[(i)] $d_0f(x):=f(p(x))$, for $f\in\linfty(G/H)$, and
\item[(ii)] $d_nf(x)=\sum_{i=1}^{n+1}(-1)^{i-1}f\big(p_{n,i}(x)\big)$, 
for $f\in\linfty\big(\fib\big)$ and $n\geq1$,
\end{enumerate}
where 
\begin{equation}\label{eq:pni}
p_{n,i}:(G/Q)^{n+1}_f\to\fib
\end{equation} 
is obtained by leaving out the $i$-th coordinate.
Observe that from the equality $(p_{n,i})_*(\nu_{n+1})=\nu_n$, 
it follows that $d_n$ is a well defined linear map
between $\linfty$ spaces. 

Then:

\begin{prop}\label{prop:resol} Let $L\leq G$ be a closed subgroup.
\begin{enumerate}
\item The complex 
\begin{equation*}
\xymatrix@1{
 0\ar[r]
&\linfty(G/H)\ar[r]
&\dots\ar[r]
&\linfty\big(\fib\big)\ar[r]^-{d_n}
&\linfty\big((G/Q)^{n+1}_f\big)\ar[r]
&\dots
}
\end{equation*}
is a strong resolution of the coefficient $L$-module $\linfty(G/H)$
by Banach $L$-modules.
\item If $Q$ is amenable and $n\geq1$, 
then the $G$-action on $\fib$ is amenable 
and $\linfty\big(\fib\big)$ is a relatively injective Banach $L$-module.
\end{enumerate}
\end{prop}

Using \cite[Theorem~2]{Burger_Monod_GAFA} (see also \S~\ref{sec:prelim}), 
this implies immediately the following:

\begin{cor}\label{cor:resol} 
Assume that $Q$ is amenable. Then the cohomology of the
complex of $L$-invariants
\begin{equation*}
\xymatrix{
 0\ar[r]
&\linfty(G/Q)^L\ar[r]
&\linfty\big((G/Q)^2_f\big)^L\ar[r]
&\dots
}
\end{equation*}
is canonically isomorphic to the bounded
continuous cohomology $\hcb^\bullet\big(L,\linfty(G/H)\big)$ 
of $L$ with coefficients in $\linfty(G/H)$.
\end{cor}

\begin{rem}\label{rem:alt} Just like for the usual resolutions of $\linfty$
functions on the Cartesian product of copies of an amenable
space (see \S~\ref{sec:prelim} or \cite{Burger_Monod_GAFA}), 
it is easy to see that 
the statements of Proposition~\ref{prop:resol} and of
Corollary~\ref{cor:resol} hold verbatim if we consider
instead the complex $\big(\la\big((G/Q)^\bullet\big),d^\bullet\big)$,
where $\la\big(\fib\big)$ is the subspace consisting of functions
in $\linfty\big(\fib\big)$ which are alternating
(observe that the symmetric group in $n$ letters acts
on $\fib$).
\end{rem}

\begin{proof}[Proof of Proposition~\ref{prop:resol}(2)]
If $n\geq1$ we have by definition the inclusion $\fib\subset(G/Q)^n$ and hence
there is a map of $G$-spaces
\begin{equation*}
\pi:\fib\to G/Q,
\end{equation*}
obtained by projection of the first component.
Since $\pi_\ast(\nu_n)=\nu$, $\pi$ realizes the measure $G$-space $\fib$ 
as an extension of the measure $G$-space $G/Q$.  If $Q$ is amenable,
the latter is an amenable $G$-space and hence the $G$-space $\fib$ 
is also amenable \cite{Zimmer:Poisson}.  Since $L$ is a closed
subgroup, $\fib$ is also an amenable $L$-space \cite[Theorem~4.3.5]{Zimmer_book} 
and hence $\linfty\big(\fib\big)$ is a relatively injective $L$-module,
\cite{Burger_Monod_GAFA}.
\end{proof}

The proof of Proposition~\ref{prop:resol}(1) consists in the
construction of appropriate contracting homotopy
operators.  Since it is rather long and technical, 
it will be given in the appendix at the end of this paper.

\subsubsection{An Implementation of the Transfer Map}\label{subsec:another_T_b} 
The point of this subsection is to see how the transfer map
in \cite{Monod_book} can be implemented, in a certain sense, 
on the resolution by $\linfty$ functions on the fibered product defined
in \S~\ref{subsec:fibered}.  

\bigskip

Let $G$ be a locally compact second countable group and $L<G$ a closed subgroup.  
The injection $L\hookrightarrow G$
gives by contravariance the restriction map
\begin{equation*}
\mathrm{r}^\bullet_\RR:\hcb^\bullet(G,\RR)\to\hcb^\bullet(L,\RR)
\end{equation*}
in bounded cohomology.
If we assume that $L\backslash G$ has a $G$-invariant probability measure $\mu$,
then the {\it transfer map}
\begin{equation*}
\mathrm{T}^\bullet:\cb(G^\bullet)^L\to\cb(G^\bullet)^G\,,
\end{equation*}
defined by integration
\begin{equation}\label{eq:transfer}
\mathrm{T}^{(n)}f(g_1,\dots,g_n):=\int_{L\backslash G}f(gg_1,\dots,gg_n)d\mu(\dot g)\,,
\end{equation}
for all $(g_1,\dots,g_n)\in G^n$,
induces in cohomology a left inverse of $\mathrm{r}^\bullet_\RR$ 
of norm one
\begin{equation*}
\mathrm{T}_\mathrm{b}^\bullet:\hcb^\bullet(L,\RR)\to\hcb^\bullet(G,\RR)\,,
\end{equation*}
(see \cite[Proposition~8.6.2,~pp.106-107]{Monod_book}).

Notice that, when dealing with the transfer map, 
the functorial machinery recalled in \S~\ref{sec:prelim} does not apply directly,
because $\mathrm{T}^\bullet$ is not a map of resolutions but 
is defined only on the subcomplex of invariant vectors. 

Let now $H,Q$ be closed subgroups of $G$ such that $Q<H$.  
We assume that $Q$ is amenable so that, by Proposition~\ref{prop:resol},
the complex $\big(\linfty\big((G/Q)^\bullet_f\big),d^\bullet\big)$ 
is a strong resolution
of the coefficient module $\linfty(G/H)$ 
by relatively injective $L$-modules. 
For $n\geq1$, $f\in\linfty\big(\fib\big)^L$, 
and $(x_1,\dots,x_n)\in\fib$, let
\begin{equation}\label{eq:tau-on-G/Q}
  (\tau^{(n)}_{G/Q}f)(x_1,\dots,x_n)
:=\int_{L\backslash G}f(gx_1,\dots,gx_n)d\mu(\dot g)\,.
\end{equation}
This defines a morphism of complexes
\begin{equation*}
\tau^\bullet_{G/Q}:\big(\linfty\big((G/Q)^\bullet_f\big)^L\big)
                 \to\big(\linfty\big((G/Q)^\bullet_f\big)^G\big)
\end{equation*}
and gives a left inverse to the inclusion
\begin{equation*}
\big(\linfty\big((G/Q)^\bullet_f\big)^G\big)\hookrightarrow\big(\linfty\big((G/Q)^\bullet_f\big)^L\big)\,.
\end{equation*}
The induced map in cohomology
\begin{equation*}
\tau_{G/Q}^\bullet:\hb^\bullet\big(L,\linfty(G/H)\big)\to
\hcb^\bullet\big(G,\linfty(G/H)\big)
\end{equation*}
is thus a left inverse of the restriction map 
$\mathrm{r}^\bullet_{\linfty(G/H)}$.

\begin{lemma}\label{lem:comm}
With the above notations, and for any amenable group $Q$, the diagram
\begin{equation}\label{eq:diagr_transfer}
\xymatrix{
 \hcb^\bullet(L,\RR)\ar[r]^{\mathrm{T}_\mathrm{b}^\bullet}\ar[d]_{\theta^\bullet_L}
&\hcb^\bullet(G,\RR)\ar[d]^{\theta^\bullet_G}\\
 \hcb^\bullet\big(L,\linfty(G/H)\big)\ar[r]^{\tau_{G/Q}^\bullet}
&\hcb^\bullet\big(G,\linfty(G/H)\big)
}
\end{equation}
commutes, where %$\mathrm{T}_\mathrm{b}^\bullet$ is the transfer map
%defined in (\ref{eq:transfer}) and 
$\theta^\bullet$ is the canonical
map induced in cohomology by the morphism of 
coefficients $\theta:\RR\to\linfty(G/H)$.
\end{lemma}

Observe that if in the above lemma we take $H=G$,
then the fibered product $\fib$ becomes the usual Cartesian
product $(G/Q)^n$, and the cohomology of the complex
of $L$-invariant $\big(\linfty\big((G/Q)^\bullet\big)^L,d^\bullet\big)$ 
computes as usual the bounded cohomology of $L$ with trivial coefficients.  
Hence we can record the following particular case of
Lemma~\ref{lem:comm}:

\begin{lemma}\label{lem:comm2}
With the above notations and for any amenable subgroup $Q\leq G$, let 
\begin{equation}\label{eq:tgq}
\mathrm{T}^\bullet_{G/Q}:
  \big(\linfty\big((G/Q)^\bullet\big)^L,d^\bullet\big)\to
  \big(\linfty\big((G/Q)^\bullet\big)^G,d^\bullet\big)
\end{equation}
be defined by 
\begin{equation}\label{eq:T_G/Q}
\mathrm{T}_{G/Q}^{(n)}f(x_1,\dots,x_n):=
\int_{L\backslash G} f(gx_1,\dots,gx_n)d\mu(g)\,,
\end{equation}
for $(x_1,\dots,x_n)\in\fib$.
Then the diagram 
\begin{equation}\label{eq:tgp-diag}
\xymatrix{
 \hcb^\bullet(L,\RR)\ar[r]^{\mathrm{T}_\mathrm{b}^\bullet}\ar[d]_{\cong}
&\hcb^\bullet(G,\RR)\ar[d]^{\cong}\\
 \hcb^\bullet(L,\RR)\ar[r]^{\mathrm{T}_{G/Q}^\bullet}
&\hcb^\bullet(G,\RR)
}
\end{equation}
commutes, where the vertical arrows are the canonical isomorphisms
in bounded cohomology extending the identity $\RR\to\RR$.
\end{lemma}

%\begin{rem} One word of comment is in order.  
%Lemma~\ref{lem:comm2} holds in full generality if $G/Q$
%is replaced by any regular amenable $G$-space.  It is natural
%to ask then why proving Lemma~\ref{lem:comm} rather than
%a similar result asserting the commutativity of a diagram
%as in (\ref{eq:diagr_transfer}) in which we consider
%general coefficients and complexes on amenable spaces.
%The reason lies in the fact that the subtlety of the situation 
%would not be captured by usual methods, and we would have needed
%anyway an {\it ad hoc} method to deal with the particular complex 
%of $\linfty$ functions on fibered products, unless we had developed
%a more general method considering complexes of $\linfty$-functions
%on appropriate sequences $(\Ss_n,\nu_n)$ of amenable spaces 
%introducing analogues of simplicial sets in the category
%of measured spaces (see \cite{Dupont}).  This would have taken 
%us however a little too far away, so we shall prove
%what is sufficient for our needs.
%\end{rem}

\begin{proof}[Proof of Lemma~\ref{lem:comm}] 
Let $G^n_f$ be the $n$-fold fibered product with respect to the
projection $G\to G/H$.  The restriction of continuous functions
defined on $G^n$ to the subspace $G^n_f\subset G^n$ induces
a morphism of strong $L$-resolutions by $L$-injective modules
\begin{equation*}
R^\bullet:\cb(G^\bullet)\to\linfty(G^\bullet_f)
\end{equation*}
extending $\theta:\RR\to\linfty(G/H)$,
so that the diagram
\begin{equation}\label{eq:d1}
\xymatrix{
\cb(G^n)^L\ar[r]^{\mathrm{T}^{(n)}}\ar[d]_{R^{(n)}_L}
&\cb(G^n)^G\ar[d]^{R^{(n)}_G}\\
 \linfty(G^n_f)^L\ar[r]^{\tau_G^{(n)}}
&\linfty(G^n_f)^G}
\end{equation} 
commutes.

Likewise, the projection $\beta_n:G^n_f\to\fib$, for $n\geq1$, gives 
by precomposition a morphism
of strong $L$-resolutions by $L$-injective modules
\begin{equation*}
\beta^\bullet:\linfty((G/Q)^\bullet_f)\to\linfty(G^\bullet_f)
\end{equation*}
extending the identity $\linfty(G/H)\to\linfty(G/H)$ and,
as before, the diagram
\begin{equation}\label{eq:d2}
\xymatrix{
 \linfty(G^n_f)^L\ar[r]^{\tau_G^{(n)}}
&\linfty(G^n_f)^G\\
 \linfty\big((G/Q)^n_f\big)^L\ar[u]^{\beta^{(n)}_L}\ar[r]^{\tau_{G/Q}^{(n)}}
&\linfty\big((G/Q)^n_f\big)^G\ar[u]_{\beta^{(n)}_G}\,,
}
\end{equation}
commutes.

The composition of the map induced in cohomology by $R^\bullet$
with the inverse of the isomorphism 
induced by $\beta^\bullet$ in cohomology realizes therefore the canonical map
\begin{equation}\label{eq:alpha}
\theta^\bullet_L:\hcb^\bullet(L,\RR)\to\hcb^\bullet\big(L,\linfty(G/H)\big)
\end{equation}
induced by the change of coefficient 
$\theta:\RR\to\linfty(G/H)$, \cite[Proposition~8.1.1]{Monod_book}.
Hence the commutative diagrams induced in cohomology by (\ref{eq:d1}) and (\ref{eq:d2})
can be combined to obtain a diagram
\begin{equation*}
\xymatrix{
 \hcb^\bullet(L,\RR)\ar[r]^{\mathrm{T}_\mathrm{b}^\bullet}\ar[d]_{R^\bullet_L}
   \ar@/_5pc/[dd]_{\theta^\bullet_L}
&\hcb^\bullet(G,\RR)\ar[d]^{R^\bullet_G}\ar@/^5pc/[dd]^{\theta^\bullet_G}\\
 \hcb^\bullet\big(L,\linfty(G/H)\big)\ar[r]^{\tau_G^\bullet}   
   \ar[d]_{(\beta^\bullet_L)^{-1}}^\cong
&\hcb^\bullet\big(G,\linfty(G/H)\big)\ar[d]^{(\beta^\bullet_G)^{-1}}_\cong\\
 \hcb^\bullet\big(L,\linfty(G/H)\big)\ar[r]^{\tau_{G/Q}^\bullet}
&\hcb^\bullet\big(G,\linfty(G/H)\big)
}
\end{equation*}
whose commutativity completes the proof.
\end{proof}

\subsubsection{An Implementation of the Pullback}\label{subsec:gpullback}
In this section we shall use the results of \cite{Burger_Iozzi_app}
to implement the pullback in bounded cohomology followed by the change
of coefficients, by using the resolution by 
$\linfty$ functions on the fibered product.  

\bigskip

Let $G'$ be a locally compact second countable group acting
on a measurable space $X$.  It is shown in \cite[Proposition~2.1]{Burger_Iozzi_app} 
that the complex $\Bb^\infty(X^\bullet)$ of bounded measurable functions
is a strong resolution of $\RR$.  Not knowing whether the modules
are relatively injective, we cannot conclude that the cohomology
of this complex computes the bounded continuous cohomology of $G'$,
however we can deduce the existence of a functorially defined map
\begin{equation*}
\epsilon^\bullet_X:\h^\bullet\big(\Bb^\infty(X^\bullet)^{G'}\big)\to\hcb^\bullet(G',\RR)
\end{equation*}
such that to any bounded measurable $G'$-invariant cocycle $c:X^{n+1}\to\RR$
corresponds canonically a class $[c]\in\hcb^n(G',\RR)$,
\cite[Corollary~2.2]{Burger_Iozzi_app}.

Let now $G$ be a locally compact second countable group, $L\leq G$
a closed subgroup acting measurably on $X$ via a continuous
homomorphism $\rho:L\to G'$,
and let us assume that there exists an $L$-equivariant measurable map
$\varphi:G/P\to X$, where $P<G$ is a closed subgroup.
The main point of \cite{Burger_Iozzi_app} is to show that the map
$\varphi$ can be used to implement the composition
\begin{equation*}
\xymatrix@1{
 \h^\bullet\big(\Bb^\infty(X^\bullet)^{G'}\big)\ar[r]^-{\epsilon^\bullet_X}
&\hcb^\bullet(G',\RR)\ar[r]^-{\rho^\bullet_\mathrm{b}}
&\hcb^\bullet(L,\RR)\,.
}
\end{equation*}
More specifically, we recall here for later use
that if $\kappa\in\hcb^n(G',\RR)$ is representable
by a $G$-invariant bounded strict measurable cocycle
$c\in\Bb^\infty(X^{n+1})^{G'}$, then the image of the pullback
$\rho_\mathrm{b}^{(n)}(\kappa)\in\hcb^n(L,\RR)$ can be represented
canonically by the cocycle in $\linfty\big((G/P)^{n+1}\big)^L$
defined by 
\begin{equation}\label{eq:old-rep}
(x_0,\dots,x_n)\mapsto c\big(\varphi(x_0),\dots,\varphi(x_n)\big)\,.
\end{equation}

The point of this section is to move one step further and
to show how to represent canonically the composition of the
above maps with the map $\theta^\bullet_L$ in (\ref{eq:alpha}).

To this purpose, let $Q,H,P$ be closed subgroups of $G$ such that 
$Q\leq H\cap P$, and let us consider the map 
\begin{equation*}
\begin{aligned}
G\times H/Q&\to G/P\\
(g,xQ)\,\,\,&\mapsto gxP
\end{aligned}
\end{equation*}
which, composed with $\varphi$, gives a measurable map
$\widetilde\varphi:G\times H/Q\to X$ which has the properties of being:
\begin{enumerate}
\item[(i)] $L$-equivariant with respect to the action 
by left translations on the first variable:
$\widetilde\varphi(\gamma g,\dot x)=\rho(\gamma)\widetilde\varphi(g,\dot x)$
for all $\gamma\in L$ and a.\,e. $(g,\dot x)\in G\times H/Q$;
\item[(ii)] $H$-invariant with respect to the right action 
considered in (\ref{eq:Haction}):
$\widetilde\varphi(gh^{-1},h\dot x)=\widetilde\varphi(g,\dot x)$
for all $h\in H$ and all $(g,\dot x)\in G\times H/Q$.
\end{enumerate}
For every $n\geq1$, the measurable map 
\begin{equation*}
\begin{aligned}
\widetilde\varphi^n_f:G\times(H/Q)^n&\longrightarrow\quad\qquad X^n\\
(g,\dot x_1,\dots,\dot x_n)&\mapsto\big(\widetilde\varphi(g,\dot x_1),\dots,\widetilde\varphi(g,\dot x_n)\big)
\end{aligned}
\end{equation*}
gives, in view of (\ref{eq:homeo-fib}), (i) and (ii), 
a measurable $L$-equivariant map $\varphi_f^n:\fib\to X^n$
defined by the composition
\begin{equation}\label{eq:varphifn}
\xymatrix@1{
 \varphi_f^n:\fib\ar[r]^-{q_n^{-1}}
&\big(G\times(H/Q)^n\big)/H\ar[r]^-{\widetilde\varphi_f^n}
&X^n\,,
}
\end{equation}
such that for every $1\leq i\leq n+1$ the diagram
\begin{equation*}
\xymatrix{
 (G/Q)^{n+1}_f\ar[r]^-{\varphi^{n+1}_f}\ar[d]_{p_{n,i}}
&X^{n+1}\ar[d]\\
 \fib\ar[r]^{\varphi^n_f}
&X^n
}
\end{equation*}
commutes, where $p_{n,i}$ was defined in (\ref{eq:pni}) and
the second vertical arrow is the map
obtained by dropping the $i$-th coordinate.
Precomposition by $\varphi^n_f$ gives thus rise to a morphism
of strong $L$-resolutions
\begin{equation*}
\xymatrix{
 0\ar[r]
&\RR\ar@{^(->}[d]\ar[r]
&\dots\ar[r]
& \Bb^\infty(X^n)\ar[d]^-{\varphi^{(n)}_f}\ar[r]
&\dots\\
 0\ar[r]
&\linfty(G/H)\ar[r]
&\dots\ar[r]
&\linfty(\fib)\ar[r]
&\dots
}
\end{equation*}
extending the inclusion $\RR\hookrightarrow\linfty(G/H)$.
Let us denote by 
\begin{equation}\label{eq:varphifbullet}
\varphi_f^\bullet:\h^\bullet\big(\Bb^\infty(X^\bullet)^{G'}\big)\to\hcb^\bullet\big(L,\linfty(G/H)\big)
\end{equation}
the map obtained in cohomology.

\begin{prop}\label{lem:4.1} Assume that $Q$ is amenable.  
Then the map $\varphi_f^\bullet$ defined in (\ref{eq:varphifbullet})
coincides with the composition
\begin{equation*}
\xymatrix@1{
 \h^\bullet\big(\Bb^\infty(X^\bullet)^{G'}\big)\ar[r]^-{\epsilon^\bullet_X}
&\hcb^\bullet(G',\RR)\ar[r]^-{\rho_\mathrm{b}^\bullet}
&\hcb^\bullet(L,\RR)\ar[r]^-{\theta^\bullet_L}
&\hcb^\bullet\big(L,\linfty(G/H)\big)\,.
}
\end{equation*}
\end{prop}

\begin{proof} By Proposition~\ref{prop:resol} 
$\big(\linfty\big((G/Q)^\bullet_f\big),d^\bullet\big)$  
is a strong resolution by relatively injective $L$-modules,
so it is enough to apply
Lemma~\ref{lem:coeff-gp} with $G=L$, $E=\linfty(G/H)$, 
$F=\RR$ the trivial coefficient $G'$-module, $F_\bullet=\Bb^\infty(X^\bullet)$,
and $E_\bullet=\big(\linfty(G/Q)^\bullet_f\big)$.
\end{proof}

For further use we record the explicit reformulation of the above
proposition:

\begin{cor}\label{cor:representative}
Let $G,G'$ be locally compact second countable groups,
$L,H,Q,P\leq G$ closed subgroups with $Q\leq H\cap P$, 
and assume that $Q$ is amenable.
Let $\rho:L\to G'$ be a continuous homomorphism, $X$ a measurable $G'$-space
and assume that there is an $L$-equivariant measurable map
$\varphi:G/P\to X$.  
Let $\kappa'\in\hcb^n(G',\RR)$
be a bounded cohomology class which admits as representative
a bounded strict $G'$-invariant measurable cocycle $c':X^{n+1}\to\RR$.
Then the class 
\begin{equation*}
\theta^{(n)}_L\big(\rho^{(n)}_\mathrm{b}(\kappa')\big)\in\hcb^n\big(L,\linfty(G/H)\big)
\end{equation*}
is represented by the $L$-invariant essentially bounded measurable cocycle
\begin{equation*}
\tilde{c}':(G/Q)_f^{n+1}\to\RR
\end{equation*}
defined by 
\begin{equation}\label{eq:ctilde}
\tilde{c}'(x_0,x_1,\dots,x_n):=c'\big(\varphi_f^n(x_0,x_1,\dots,x_n)\big)\,,
\end{equation}
where $\varphi_f^n$ is defined in (\ref{eq:varphifn}).

In particular, if $L=G=G'$, $\rho=\id$ and $X=G/P$, then
if the class $\kappa\in\hcb^n(G,\RR)$ admits as representative 
a bounded strict $G$-invariant Borel cocycle
$c:(G/P)^{n+1}\to\RR$, the class
\begin{equation*}
\theta^{(n)}_G(\kappa)\in\hcb^{(n)}\big(G,\linfty(G/H)\big)\,.
\end{equation*}
is represented by the bounded strict $G$-invariant Borel cocycle
\begin{equation*}
\tilde c:(G/Q)_f^{n+1}\to\RR
\end{equation*}
defined by 
\begin{equation}\label{eq:ctilde-on-G/P}
\tilde c(x_1,\dots,x_{n+1}):=c(x_1P,\dots,x_{n+1}P)\,.
\end{equation}
\end{cor}

%%% Local Variables: 
%%% mode: latex
%%% TeX-master: "deform"
%%% End: 

\vskip1cm
\subsection{A Formula}\label{subsec:formula} 
We apply now all the results obtained so far 
to prove the first instance of a formula deriving
from a commutative diagram in cohomology.  The main
application of this formula will involve a numerical invariant 
attached to a representation, but in order to obtain this 
we need to know that some of the cohomology spaces involved 
are one-dimensional, like for instance if 
the groups $G$ and $G'$ are of Hermitian type and simple
(see \S~\ref{subsec:hermitian}).  We want however 
to isolate here the preliminary result which is true
for all locally compact second countable groups 
and which could be of independent interest.

\bigskip

\begin{prop}\label{prop:formula-homog}  Let $G$ and $G'$ be locally compact
second countable groups, let 
$L,H,Q,P\leq G$ be closed subgroups with $Q\leq H\cap P$, 
and assume $L\backslash G$ carries a $G$-invariant probability measure
$\mu$.  
Let $\rho:L\to G'$ be a continuous homomorphism, $X$ a measurable $G'$-space
and assume that there is an $L$-equivariant measurable map
$\varphi:G/P\to X$.  Let $\kappa'\in\hcb^2(G',\RR)$ and let 
$\kappa:=\mathrm{T}_\mathrm{b}^{(2)}\big(\rho^{(2)}_\mathrm{b}(\kappa')\big)
\in\hcb^2(G,\RR)$.
Let $c\in\Bb^\infty\big((G/P)^3\big)^{G}$ and 
$c'\in\Bb^\infty(X^3)^{G'}$ be alternating cocycles representing 
$\kappa$ and $\kappa'$ respectively, and
let $\tilde{c}:(G/Q)^3_f\to\RR$ and $\tilde{c}':(G/Q)^3_f\to\RR$
be the corresponding alternating cocycles defined respectively 
in (\ref{eq:ctilde-on-G/P}) and (\ref{eq:ctilde}).

Assume that $Q$ is amenable and that $H$ acts ergodically on $H/Q\times H/Q$.  
Then, if $\varphi^3_f$ is the map defined in 
(\ref{eq:varphifn}), we have
\begin{equation*}
\int_{L\backslash G}\tilde{c}'\big(\varphi^3_f(gx_1,gx_2,gx_3)\big)d\mu(\dot g)=
\tilde{c}(x_1,x_2,x_3)
\end{equation*}
for a.\,e. $(x_1,x_2,x_3)\in(G/Q)^3_f$. 
\end{prop}

\begin{proof} The commutativity of the square in the following 
diagram (see Lemma~\ref{lem:comm})
\begin{equation*}
\xymatrix{
 \h^2\big(\Bb^\infty(X^3)^{G'}\big)\ar[r]^-{\omega^{(2)}_X}
&\hcb^2(G',\RR)\ar[r]^-{\rho_\mathrm{b}^{(2)}}
&\hcb^2(L,\RR)\ar[r]^-{\theta^{(2)}_L}\ar[d]^{\mathrm{T}_\mathrm{b}^{(2)}}
&\hcb^2\big(L,\linfty(G/H)\big)\ar[d]^{\tau_{G/Q}^{(2)}}\\
\h^2\big(\Bb^\infty((G/P)^3)^{G}\big)\ar[rr]^-{\omega^{(2)}_{G/P}}
&
&\hcb^2(G,\RR)\ar[r]^-{\theta^{(2)}_G}
&\hcb^2\big(G,\linfty(G/H)\big)
}
\end{equation*}
applied to the class $\rho_\mathrm{b}^{(2)}(\kappa')\in\hcb^2(L,\RR)$
reads
\begin{equation*}
 \tau_{G/Q}^{(2)}\big(\theta_L^{(2)}\big(\rho_\mathrm{b}^{(2)}(\kappa')\big)\big)
=\theta_G^{(2)}\big(\mathrm{T}_\mathrm{b}^{(2)}\big(\rho_\mathrm{b}^{(2)}(\kappa')\big)\big)
=\theta_G^{(2)}(\kappa)\,.
\end{equation*}
Hence the representatives for the classes $\theta_G^{(2)}(\kappa)$ and
$\theta_L^{(2)}\big(\rho_\mathrm{b}^{(2)}(\kappa')\big)$
chosen according to 
Corollary~\ref{cor:representative} satisfy the relation
\begin{equation*}
\tau_{G/Q}^{(2)}(\tilde{c}')=\tilde{c}+db\,,
\end{equation*}
where $b\in\linfty\big((G/Q)^2_f\big)^G$,
which, using the definition of $\tau_{G/Q}^{(2)}$
in (\ref{eq:tau-on-G/Q}) implies that 
\begin{equation*}
\int_{L\backslash G}\tilde{c}'\big(\varphi^3_f(gx_1,gx_2,gx_3)\big)d\mu(\dot g)=
\tilde{c}(x_1,x_2,x_3)+db
\end{equation*}
for a.\,e. $(x_1,x_2,x_3)\in(G/Q)^3_f$. 
However $G$ acts ergodically on $(G/Q)^2_f$ because
it acts on the basis of the fibration $(G/Q)^2_f\to G/H$ 
transitively with stabilizer $H$, which then 
by hypothesis acts ergodically on the typical fiber 
homeomorphic to $(H/Q)^2$.  Hence $\linfty((G/Q)^2_f)^G=\RR$.
Thus $db$ is constant and hence zero, being the
difference of two alternating functions.
\end{proof}

Notice that if $G=H$ and $P=Q$, the above
formula takes a more familiar form, namely

\begin{cor}\label{cor:usual-formula} Let $\kappa'\in\hcb^2(G',\RR)$
and $\kappa:=\mathrm{T}_\mathrm{b}^{(2)}\big(\rho^{(2)}_\mathrm{b}(\kappa')\big)$.
With the same notation as in Corollary~\ref{cor:representative}
assume that $L\backslash G$ carries a $G$-invariant probability
measure $\mu$.  Let $c\in\Bb_\mathrm{alt}^\infty\big((G/P)^3\big)^{G}$ and
$c'\in\Bb_\mathrm{alt}^\infty(X^3)^{G'}$ represent $\kappa$ and $\kappa'$
respectively and 
let $\tilde{c}:(G/P)^3\to\RR$ and $\tilde{c}':(G/P)^3\to\RR$
be the corresponding alternating cocycles defined in (\ref{eq:ctilde-on-G/P}) and 
(\ref{eq:ctilde}).

Assume that $G$ acts ergodically on $G/P\times G/P$.  Then
\begin{equation*}
\int_{L\backslash G}
  c'\big(\varphi(x_1),\varphi(x_2),\varphi(x_3)\big)d\mu(\dot g)=
  c(x_1,x_2,x_3)
\end{equation*}
for a.\,e. $(x_1,x_2,x_3)\in(G/P)^3$.
\end{cor}
%%% Local Variables: 
%%% mode: latex
%%% TeX-master: "deform"
%%% End: 

\vskip1cm
\subsection{The Toledo Map and the Bounded Toledo Map}\label{sec:Toledo}
We shall now see to which extent one can define a transfer map 
also in ordinary continuous cohomology and link it to the one previously
defined in bounded cohomology.  As we shall see, problems arise 
if the subgroup $L$ is only of finite covolume and not cocompact.

\bigskip

Let $L$  be a closed subgroup of a locally compact group $G$.
Then we have a commutative diagram 
\begin{equation}\label{eq:square}
\xymatrix{
 \hcb^\bullet(L,\RR)\ar[d]_-{\imath_\mathrm{b}^\bullet}\ar[r]^{{\tt c}^\bullet_L}
& \hc^\bullet(L,\RR)\ar[d]^-{\imath^\bullet}\\
\hcb^\bullet\big(G,\linfty(L\backslash G)\big)\ar[r]^{{\tt c}^\bullet_G}
&\hc^\bullet\big(G,\mathrm{L}_{\mathrm{loc}}^p(L\backslash G)\big)
}
\end{equation}
for $1\leq p<\infty$, where $\imath^\bullet$ is the induction in ordinary
continuous cohomology (see \cite{Blanc}),
$\imath^\bullet_\mathrm{b}$ is the isometric isomorphism
defined by the induction in bounded continuous cohomology
(see \cite{Burger_Monod_GAFA})
and the horizontal arrows are comparison maps.
(Observe that for trivial coefficients the comparison map 
was recalled in the introduction; for a thorough discussion of 
this diagram see \cite[\S~10.1]{Monod_book}.)
If $L\backslash G$ carries a $G$-invariant probability measure $\mu$,
we have, for $1\leq p<\infty$, an obvious morphism of coefficient modules
\begin{equation*}
\begin{aligned}
\mathrm{m}:\mathrm{L}^p(L\backslash G)&\to\quad\RR\\
f\qquad&\mapsto\int_{L\backslash G} fd\mu\,.
\end{aligned}
\end{equation*}
If moreover $L\backslash G$ is compact, then 
$\mathrm{L}_\mathrm{loc}^p(L\backslash G)=\mathrm{L}^p(L\backslash G)$,
and composing $\imath$ with the change of coefficients $\mathrm{m}$ in continuous 
cohomology gives a transfer map in ordinary continuous cohomology
\begin{equation*}\label{eq:transfer-ord}
\mathrm{T^\bullet}:\hc^\bullet(L,\RR)\to\hc^\bullet(G,\RR)
\end{equation*}
which is a left inverse to the restriction map and 
leads to a commutative diagram
\begin{equation}\label{eq:transfer-bounded-ordinary}
\xymatrix{
 \hcb^\bullet(L,\RR)\ar[d]_{\mathrm{T}_\mathrm{b}^\bullet}\ar[r]^{{\tt c}^\bullet_L}
&\hc^\bullet(L,\RR)\ar[d]^{\mathrm{T}^\bullet}\\
 \hcb^\bullet(G,\RR)\ar[r]^{{\tt c}^\bullet_G}
&\hc^\bullet(G,\RR)
}
\end{equation}
which is very useful in applications when it comes to identifying invariants
in bounded cohomology in terms of ordinary cohomological invariants.

\begin{rem} The above fails if $L\backslash G$ is only of finite volume.
For example, if $L=\Gamma<G$ is a nonuniform lattice,
then there is in general no left inverse to the restriction 
in cohomology $\hc^\bullet(G,\RR)\to\h^\bullet(\Gamma,\RR)$ 
as this map is often not injective.
In fact, one can for example consider the case in which
$\Xx=G/K$ is an $n$-dimensional symmetric space of noncompact type:
then $\hc^n(G,\RR)=\Omega^n(\Xx)^G$ is generated by the volume form
and hence not zero, while if $\Gamma<G$ is any nonuniform torsionfree
lattice, the cohomology $\h^n(\Gamma,\RR)$ vanishes as it is isomorphic 
to $\hdr^n(\Gamma\backslash \Xx)$.

The extent to which one can remedy this mishap by keeping track
of the relation with bounded continuous cohomology is thus the object of 
\S~\ref{sec:comp} in the context when $G$ is a semisimple Lie group.
\end{rem}

\bigskip

We apply here the above considerations to associate to a homomorphism 
a map in cohomology which will in some cases produce in addition
a numerical invariant (see for instance the Hermitian case, 
\S~\ref{subsec:hermitian}).

\bigskip

Let $L\leq G$ be a closed subgroup of a locally compact group $G$, 
and $\rho:L\to G'$ a continuous homomorphism into a locally compact group $G'$.  
The composition of the pullback
\begin{equation*}
\rho^\bullet_\mathrm{b}:\hcb^\bullet(G',\RR)\to\hcb^\bullet(L,\RR)
\end{equation*}
with the transfer map $\mathrm{T}_\mathrm{b}^\bullet$
defined in \S~\ref{subsec:another_T_b}
gives rise to the {\it bounded Toledo map}
\begin{equation}\label{eq:bounded-toledo-map}
\mathrm{T}_\mathrm{b}^\bullet(\rho):\hcb^\bullet(G',\RR)\to\hcb^\bullet(G,\RR)
\end{equation}
which provides a basic invariant of the homomorphism $\rho:L\to G'$. 

We remark once again that the bounded Toledo map is
defined for {\it all} closed subgroups of $G$
such that on $L\backslash G$ there is a $G$-invariant probability measure.
If however $L\backslash G$ is in addition  compact 
(for example, a uniform lattice)
then we also have an analogous construction in ordinary
cohomology.  Namely, associated to the homomorphism
$\rho:L\to G'$ we have a morphism
\begin{equation*}
\rho^\bullet:\hc^\bullet(G',\RR)\to\hc^\bullet(L,\RR)
\end{equation*}
which, composed with the transfer map $\mathrm{T}^\bullet$
defined in (\ref{eq:transfer-ord}) gives a map
\begin{equation*}
\mathrm{T}^\bullet(\rho):\hc^\bullet(G',\RR)\to\hc^\bullet(G,\RR)
\end{equation*}
which we denote by the ``Toledo map'' and has the property
that the diagram
\begin{equation*}%\label{eq:bdd-tol-tol}
\xymatrix{
 \hcb^\bullet(G',\RR)\ar[d]_{ \mathrm{T}_\mathrm{b}^\bullet(\rho)}\ar[r]^{{\tt c}^\bullet_{G'}}
&\hc^\bullet(G',\RR)\ar[d]^{ \mathrm{T}^\bullet(\rho)}\\
 \hcb^\bullet(G,\RR)\ar[r]^{{\tt c}^\bullet_G}
&\hc^\bullet(G,\RR)\,,
}
\end{equation*}
where the horizontal arrows are comparison maps, commutes.

The interplay between these two maps is the basic ingredient in
this paper for the cocompact case, as well as in \cite{Iozzi_ern},
\cite{Burger_Iozzi_Wienhard_ann} and \cite{Burger_Iozzi_Wienhard_tol}.
In the finite volume case we will need to resort to a somewhat 
more elaborate version of the above diagram which can be developed
when $G$ is a connected semisimple Lie group -- see (\ref{eq:big-one}).
%%% Local Variables: 
%%% mode: latex
%%% TeX-master: "deform"
%%% End: 

\vskip1cm
\section{More Cohomological Tools for Semisimple Lie Groups and Applications}
In this section we specialize the discussion to semisimple Lie groups
(with finite center) and their closed subgroups.  The immediate
advantage will be the identification of the ordinary group cohomology 
with the cohomology of differential forms on symmetric spaces.
\subsection{A Factorization of the Comparison Map}\label{sec:comp}
The main point of this section is to provide, in the case of semisimple Lie groups,
a substitute to the the missing arrow in 
\begin{equation}\label{eq:missing}
\xymatrix{
 \hcb^\bullet(L,\RR)\ar[d]_{\mathrm{T}_\mathrm{b}^\bullet}\ar[r]^{{\tt c}^\bullet_L}
&\hc^\bullet(L,\RR)\ar@{.>}[d]\\
 \hcb^\bullet(G,\RR)\ar[r]^{{\tt c}^\bullet_G}
&\hc^\bullet(G,\RR)
}
\end{equation}
if the subgroup $L\leq G$ is only of finite covolume.

\bigskip

Let $G$ be a connected semisimple Lie group with finite center
and $\Xx$ the associated symmetric space.
Any closed subgroup $L\leq G$ acts properly on $\Xx$ and hence the complex
\begin{equation*}
\xymatrix@1{
 \RR\ar[r]
&\Omega^0(\Xx)\ar[r]
&\dots\ar[r]
&\Omega^k(\Xx)\ar[r]
&\dots
}
\end{equation*}
of $\mathrm{C}^\infty$ differential forms on $\Xx$ with the usual 
exterior differential is a resolution by continuous injective
$L$-modules (where injectivity now refers to the usual notion 
in continuous cohomology), from which one obtains a canonical isomorphism
\begin{equation*}
\xymatrix@1{
\hc^\bullet(L,\RR)\ar[r]^-\cong
&\h^\bullet\big(\Omega^\bullet(\Xx)^L\big)
}
\end{equation*}
in cohomology, \cite{Mostow_61}.
Let moreover $\big(\Omega_\infty^\bullet(\Xx),d^\bullet\big)$ denote the complex
of smooth differential forms $\alpha$ on $\Xx$ such that
$x\mapsto\|\alpha_x\|$ and $x\mapsto\|d\alpha_x\|$ are in $\linfty(\Xx)$,
and let $h(\Xx)$ denote the volume entropy of $\Xx$, that is the rate of
exponential growth of volume of geodesic balls in $\Xx$, \cite{Eberlein}.
Then
\begin{prop}\label{prop:transfer-L}
Let $G$ be a connected semisimple Lie group with finite center,
$\Xx$ the associated symmetric space and $L\leq G$ any closed subgroup.
Then there exists a map
\begin{equation*}
\delta^\bullet_{\infty,L}:\hcb^\bullet(L,\RR)\to\h^\bullet\big(\Omega^\bullet_\infty(\Xx)^L\big)
\end{equation*}
such that the diagram 
\begin{equation}\label{eq:transfer-L}
\xymatrix{
 \hcb^\bullet(L,\RR)\ar[rr]^{{\tt c}^\bullet_L}\ar[rrd]_{\delta^\bullet_{\infty,L}}
&
&\hc^\bullet(L,\RR)\ar[r]^\cong
&\h^\bullet\big(\Omega^\bullet(\Xx)^L\big)\\
&
&\h^\bullet\big(\Omega^\bullet_\infty(\Xx)^L\big)\,,\ar[ur]_{i^\bullet_{\infty,L}}
&
}
\end{equation}
commutes, where $i^\bullet_{\infty,L}$ is 
the map induced in cohomology by the inclusion of complexes
\begin{equation*}
i^\bullet_\infty:\Omega_\infty^\bullet(\Xx)\to\Omega^\bullet(\Xx)\,.
\end{equation*}
Moreover, the norm of $\delta^{(k)}_{\infty,L}$ is bounded by $h(\Xx)^k$.
\end{prop}

Before proving the proposition, we want to push our result  
a little further in the case when $L=\Gamma<G$ is a lattice.  
In particular, we are going to see how the map
$\delta^\bullet_{\infty,\Gamma}$ fits into a diagram where the transfer 
appears.  If $1\leq p\leq\infty$, 
let $\Omega^n_p(\Xx)^\Gamma$ be the space of $\Gamma$-invariant
smooth differential $n$-forms on $\Xx$ such that
$x\mapsto\|\alpha_x\|$ and $x\mapsto\|d\alpha_x\|$ are in $\lp(\Gamma\backslash \Xx)$,
and consider the complex $\big(\Omega^\bullet_p(\Xx)^\Gamma,d^\bullet\big)$.
Incidentally, notice that this is a rather misleading notation
if $\Xx$ is not compact, because in this case only for $p=\infty$ 
one has that $\big(\Omega^\bullet_\infty(\Xx)^\Gamma,d^\bullet\big)$
is the subcomplex of invariants of $\big(\Omega^\bullet_\infty(\Xx),d^\bullet\big)$.
Let $\delta_{p,\Gamma}^\bullet$ be the map obtained by 
composing the map $\delta^\bullet_{\infty,\Gamma}$ 
in Proposition~\ref{prop:transfer-L} 
with the map obtained by the inclusion of complexes 
\begin{equation*}
\Omega^\bullet_\infty(\Xx)^\Gamma\to\Omega^\bullet_p(\Xx)^\Gamma\,,
\end{equation*}
namely
\begin{equation*}
\xymatrix{
 \hb^\bullet(\Gamma,\RR)\ar[r]^-{\delta^\bullet_{\infty,\Gamma}}\ar@/^2pc/[rr]^{\delta^\bullet_{p,\Gamma}}
&\h^\bullet\big(\Omega^\bullet_\infty(\Xx)^\Gamma\big)\ar[r]
&\h^\bullet\big(\Omega_p^\bullet(\Xx)^\Gamma\big)\,.
}
\end{equation*}
Also, since $\Omega(\Xx)^G\subset\Omega_\infty(\Xx)$ and 
$\Gamma\backslash \Xx$ is of finite volume, 
the restriction map
\begin{equation*}
\Omega^\bullet(\Xx)^G\to\Omega^\bullet_p(\Xx)^\Gamma
\end{equation*}
is defined and admits a left inverse $j^\bullet_p$ defined by integration
\begin{equation*}
j^\bullet_p\alpha=\int_{\Gamma\backslash G}(L_g^\ast\alpha)d\mu(\dot g)\,,
\end{equation*}
for $\alpha\in\Omega_p^\bullet(\Xx)^\Gamma$
and where $L_g$ is left translation by $g$.
The following proposition gives an interesting diagram to be 
compared with (\ref{eq:missing})

\begin{prop}\label{prop:transfer-Gamma} 
Let $G$ be a connected semisimple Lie group with finite center
and associated symmetric space $\Xx$, and let $\Gamma<G$ be a lattice.
The following diagram 
\begin{equation}\label{eq:bounded-rep}
\xymatrix{
 \hb^\bullet(\Gamma,\RR)\ar[dd]_{\mathrm{T}_\mathrm{b}^\bullet}
       \ar[rr]^{{\tt c}^\bullet_\Gamma}\ar[rrd]_{\delta_{p,\Gamma}^\bullet}
&
&\h^\bullet(\Gamma,\RR)\ar[r]^\cong
&\h^\bullet\big(\Omega^\bullet(\Xx)^\Gamma\big)\\
&
&\h^\bullet\big(\Omega^\bullet_p(\Xx)^\Gamma\big)\ar[ur]_{i^\bullet_{p,\Gamma}}\ar[dr]^{j^\bullet_p}
& \\
 \hcb^\bullet(G,\RR)\ar[rr]^{{\tt c}^\bullet_G}
&
&\hc^\bullet(G,\RR)\ar[r]^\cong
&\Omega^\bullet(\Xx)^G
}
\end{equation}
commutes for all $1\leq p\leq\infty$.
\end{prop}

%We saw in \S~\ref{sec:prelim} that if $\partial\Xx$ is 
%the geometric boundary of $\Xx$ (with an appropriate measure $\nu_0$ -- 
%see below), the bounded cohomology
%n$\hb^\bullet(\Gamma,\RR)$ is canonically isomorphic 
%to the cohomology of the complex
%$(\linfty((\partial\Xx)^{\bullet+1},\nu_0^{\bullet+1})^\Gamma,d^\bullet)$.
We start the proof by showing how to associate to an $\linfty$
function $c$ on $(\partial\Xx)^{n+1}$ a differential $n$-form
obtained by integrating, with respect to an appropriate density at infinity
and weighted by the function $c$,
the differential form obtained from the Busemann vectors
associated to $n$ points at infinity.

So, let $\partial \Xx$ be the geodesic ray boundary of $\Xx$ and
\begin{equation*}
B:\partial \Xx\times \Xx\times \Xx\to\RR
\end{equation*}
the Busemann cocycle. Fix a basepoint $0\in \Xx$ and let $K=\stab_G(0)$,
$\frakg=\frakk\oplus\frakp$ the associated Cartan decomposition,
$\fraka^+\subset\frakp$ a positive Weyl chamber and $b\in\fraka^+$
the vector predual to the sum of the positive roots associated to $\fraka^+$. 
Then $h(\Xx)=\|b\|$.  Let $\xi_b\in \partial \Xx$ be the point
at infinity determined by $b$;
let $\nu_0$ be the unique $K$-invariant probability measure
on $G\xi_b\subset\partial\Xx$.  Then
\begin{equation}\label{eq:busemann}
d(g_\ast\nu_0)(\xi)=e^{-hB_\xi(g0,0)}d\nu_0(\xi)\,.
\end{equation}
For $\xi\in \partial \Xx$, let us define a $\mathrm{C}^\infty$ map by
\begin{equation}\label{eq:exi}
\begin{aligned}
e^\xi:\Xx&\to\quad\quad\RR\\
x&\mapsto e^{-h(\Xx)B_\xi(x,0)}\,.
\end{aligned}
\end{equation}

\begin{lemma}\label{lem:explicit}
Let $G$ be a connected semisimple Lie group with finite center, 
and let $\Xx$ be its associated symmetric space with geodesic ray boundary 
$\partial\Xx$. For each $c\in\linfty((\partial \Xx)^{n+1},\nu_0^{n+1})$,
the differential form defined by
\begin{equation}\label{eq:explicit}
\omega:=\int_{(\partial \Xx)^{n+1}} 
     c(\xi_0,\dots,\xi_n)e^{\xi_0}\wedge de^{\xi_1}\wedge\dots\wedge de^{\xi_n}
     d\nu_0^{n+1}(\xi_0,\dots,\xi_n)\,.
\end{equation}
is in $\Omega^n_\infty(\Xx)$.
Moreover the resulting map
\begin{equation*}
\begin{aligned}
\delta^\bullet_\infty:\linfty\big((\partial \Xx)^{\bullet+1},\nu_0^{\bullet+1}\big)
   &\to\Omega^\bullet_\infty(\Xx)\\
c\qquad\quad&\mapsto\quad\omega
\end{aligned}
\end{equation*}
is a $G$-equivariant map of complexes, and
\begin{equation}\label{eq:ineq}
\|\delta^{(n)}_\infty\|\leq h(\Xx)^n\,.
\end{equation}
\end{lemma}

\begin{proof} For $\xi\in \partial \Xx$, let $X_\xi(x)$ be 
the unit tangent vector at $x$ pointing in the direction of $\xi$,
and let $g_x(\,\cdot\,,\,\cdot\,)$ be the Riemannian metric on $\Xx$ at $x$.
Since the gradient of the Busemann function $B_\xi(x,0)$ at $x$
is $-X_\xi(x)$ \cite{Eberlein}, we have that for $v\in(T\Xx)_x$,
$(dB_\xi)_x(v)=-g_x\big(v,X_\xi(v)\big)$. 
Then 
\begin{equation*}
(de^\xi)_x(v)=h(\Xx)g_x\big(v,X_\xi(x)\big)e^\xi(x)\,.
\end{equation*}
This implies that if $v_1,\dots,v_n$ are tangent vectors based at $x$, then
\begin{equation*}
\begin{aligned}
    &|\omega_x(v_1,\dots,v_n)|\\
&\leq h(\Xx)^n\int_{(\partial \Xx)^{n+1}}
        \big|c(\xi_0,\xi_1,\dots,\xi_n)\big|e^{\xi_0}(x)\cdot\\
&\hphantom{\leq h(\Xx)^n\int_{(\partial \Xx)^{n+1}}\big|c(\xi_0,}
             \cdot\left(\prod_{i=1}^ng_x\big(v_i,X_{\xi_i}(x)\big)e^{\xi_i}(x)\right)
                 d\nu_0^{n+1}(\xi_0,\dots,\xi_n)\\
&\leq h(\Xx)^n\int_{(\partial \Xx)^{n+1}}\|c\|_\infty e^{\xi_0}(x)
             \left(\prod_{i=1}^n\|v_i\|e^{\xi_i}(x)\right) 
                 d\nu_0^{n+1}(\xi_0,\dots,\xi_n)\\
&=h(\Xx)^n\|c\|_\infty\left(\int_{\partial \Xx} e^{\xi_0}(x)d\nu_0(\xi_0)\right)
             \prod_{i=1}^n
                  \left(\|v_i\|\int_{\partial \Xx} e^{\xi_i}(x)d\nu_0(\xi_i)\right)
\end{aligned}
\end{equation*}
But writing $x=g0$ and using that, as indicated in (\ref{eq:busemann}),
$d(g_\ast\nu_0)$ is a probability measure, we get that for all $0\leq i\leq n$
and all $x\in \Xx$
\begin{equation}\label{eq:int-one}
\int_{\partial \Xx} e^{\xi_i}(x)d\nu_0(\xi_i)=1\,,
\end{equation}
which shows that 
\begin{equation*}
|\omega_x(v_1,\dots,v_n)|\leq h(\Xx)^n\|c\|_\infty\prod_{i=1}^n\|v_i\|\,,
\end{equation*}
so that 
\begin{equation*}
\|\delta^{(n)}_\infty c\|=\sup_{x\in \Xx}\,\sup_{\|v_1\|,\dots,\|v_n\|\leq1}
                            |\omega_x(v_1,\dots,v_n)|
\leq h(\Xx)^n\|c\|_\infty\,.
\end{equation*}
This proves (\ref{eq:ineq}) and 
the fact that the image $\delta^{(n)}_\infty(c)$ is a bounded form.  
Once we shall have proven that $\delta^{(n)}_\infty(dc)=d\delta^{(n-1)}_\infty(c)$, 
it will follow automatically
that also $d\delta_\infty^{n-1}(c)$ is bounded and hence the image of $\delta^\bullet_\infty$ is in 
$\Omega^\bullet_\infty(\Xx)$.  
To this purpose, let us compute for $c\in\linfty((\partial \Xx)^n,\nu_0^n)$
\begin{equation*}
\begin{aligned}
\delta^{(n)}_\infty(dc)
=\int_{(\partial \Xx)^{n+1}} &e^{\xi_0}\wedge de^{\xi_1}\wedge\dots\wedge de^{\xi_n}\cdot\\
&\cdot\left[\sum_{i=0}^n(-1)^ic(\xi_0,\dots,\hat\xi_i,\dots,\xi_n)\right]
   d\nu_0^{n+1}(\xi_0,\dots,\xi_n)\,.
\end{aligned}
\end{equation*}
For $i\geq1$ the $i$-th term is
\begin{equation*}
\begin{aligned}
 &(-1)^i\int_{(\partial \Xx)^{n+1}} 
     e^{\xi_0}\wedge de^{\xi_1}\wedge\dots\wedge de^{\xi_n}\cdot\\
 &\hphantom{(-1)^i\int_{(\partial \Xx)^{n+1}}e^{\xi_0}\wedge de^{\xi_1}\wedge\dots}
     \cdot c(\xi_0,\dots,\hat\xi_i,\dots,\xi_n)d\nu_0^{n+1}(\xi_0,\dots,\xi_n)\\
=&-d\left(\int_{\partial \Xx} e^{\xi_i}d\nu_0(\xi_i)\right)\wedge\dots=0
\end{aligned}
\end{equation*}
since by (\ref{eq:int-one})
\begin{equation*}
d\left(\int_{\partial \Xx} e^{\xi_i}d\nu_0(\xi_i)\right)=0\,.
\end{equation*}
Thus
\begin{equation*}
\begin{aligned}
   \delta^{(n)}_\infty(dc)
&=\int_{(\partial \Xx)^{n+1}} e^{\xi_0}\wedge de^{\xi_1}\wedge\dots\wedge de^{\xi_n}\,
           c(\xi_1,\dots,\xi_n) d\nu_0^{n+1}(\xi_0,\dots,\xi_n)\\
&=\left[\int_{\partial \Xx} e^{\xi_0}d\nu_0(\xi_0)\right]
           \int_{(\partial \Xx)^n} de^{\xi_1}\wedge\dots\wedge de^{\xi_n}\cdot\\
&\hphantom{=\left[\int_{\partial \Xx} e^{\xi_0}d\nu_0(\xi_0)\right]
           \int_{(\partial \Xx)^n} de^{\xi_1}\wedge}\cdot c(\xi_1,\dots,\xi_n) d\nu_0^n(\xi_1,\dots,\xi_n)\\
&=\int_{(\partial \Xx)^n} de^{\xi_1}\wedge\dots\wedge de^{\xi_n}c(\xi_1,\dots,\xi_n)d\nu_0^n(\xi_1,\dots,\xi_n)\,.
\end{aligned}
\end{equation*}
On the other hand
\begin{equation*}
\delta^{(n-1)}_\infty(c)=\int_{(\partial \Xx)^n} e^{\xi_1}\wedge de^{\xi_2}\wedge\dots\wedge de^{\xi_n}\,
           c(\xi_1,\dots,\xi_n)d\nu_0^n(\xi_1,\dots,\xi_n)\,,
\end{equation*}
so that by definition
\begin{equation*}
\begin{aligned}
d\delta^{(n-1)}_\infty(c)&=\int_{(\partial \Xx)^n} de^{\xi_1}\wedge de^{\xi_2}\wedge\dots\wedge de^{\xi_n} \,
           c(\xi_1,\dots,\xi_n)d\nu_0^n(\xi_1,\dots,\xi_n)
          \\&=\delta^{(n)}_\infty(dc)\,.
\end{aligned}
\end{equation*}
The $G$-equivariance of $\delta^\bullet_\infty$ follows from (\ref{eq:busemann})
and the cocycle property of the Busemann function $B_\xi(x,y)$,
hence completing the proof.
\end{proof}

\begin{proof}[Proof of Proposition~\ref{prop:transfer-L}] 
This is a direct application of \cite[Proposition~9.2.3]{Monod_book}.
Indeed, since the $L$-action on $(\partial\Xx,\nu_0)$ is amenable,
we have that $(\linfty((\partial \Xx)^{\bullet+1},\nu_0^{\bullet+1}))$
is a strong resolution of $\RR$ by relatively injective $L$-modules
\cite{Burger_Monod_GAFA}; moreover, it is well known that, 
$(\Omega^\bullet(\Xx),d^\bullet)$ is a resolution of $\RR$ 
by injective continuous $L$-modules, 
where in this case injectivity is meant in ordinary cohomology
(see \cite{Mostow_61}), and $\Omega^\bullet(\Xx)$ is as usual equipped
with the $\mathrm{C}^\infty$-topology.  Finally one checks
on the formulas that the composition $i^\bullet_\infty\circ \delta^\bullet_\infty$, 
where $i^\bullet_\infty$ is the injection 
\begin{equation*}
i^\bullet_\infty:\Omega^\bullet_\infty(\Xx)\to\Omega^\bullet(\Xx)\,,
\end{equation*}
is a continuous $L$-morphism of complexes.
The hypotheses of \cite[Proposition~9.2.3]{Monod_book}
are hence verified and thus 
\begin{equation*}
i^\bullet_{\infty,L}\circ\delta^\bullet_{\infty,L}:
\h^\bullet\big(\linfty\big((\partial \Xx)^{\bullet+1},\nu_0^{\bullet+1}\big)^L\big)
\to
\h^\bullet\big(\Omega^\bullet(\Xx)^L\big)
\end{equation*}
realizes the canonical comparison map
$${\tt c}^\bullet_L:\hcb^\bullet(L,\RR)\to\hc^\bullet(L,\RR)$$.
\end{proof}

\begin{proof}[Proof of Proposition~\ref{prop:transfer-Gamma}]  
The proof of Proposition~\ref{prop:transfer-L} remains
valid verbatim for all $1\leq p\leq\infty$ to show
the commutativity of the upper diagram, so it remains to show 
only the commutativity of the lower part.
Notice moreover that since 
\begin{equation*}
i^\bullet_{p,G}:\Omega_p^\bullet(\Xx)^G\to\Omega^\bullet(\Xx)^G
\end{equation*}
is the identity, $\delta^\bullet_{p,G}$ realizes in cohomology
the canonical comparison map.
Furthermore, if $P$ is the minimal parabolic in $G$ stabilizing $\xi_b$
and we identify $(\partial\Xx,\nu_0)$ with $(G/P,\nu_0)$
as measure spaces, the commutativity of the diagram
\begin{equation*}
\xymatrix{
 \linfty\big((\partial \Xx)^{\bullet+1},\nu_0^{\bullet+1}\big)^\Gamma
           \ar[r]^-{\delta^\bullet_{p,\Gamma}}
           \ar[d]_{\mathrm{T}_{G/P}^{\bullet+1}}
&\Omega_p^\bullet(\Xx)^\Gamma\ar[d]^{j_p^\bullet}\\
 \linfty\big((\partial \Xx)^{\bullet+1},\nu_0^{\bullet+1}\big)^G\ar[r]^-{\delta^\bullet_{p,G}}
&\Omega^\bullet(\Xx)^G
}\,,
\end{equation*}
is immediate, where $\mathrm{T}^\bullet_{G/P}$ is defined in (\ref{eq:T_G/Q}).
The commutativity of the diagram (\ref{eq:tgp-diag}) in Lemma~\ref{lem:comm}
then completes the proof.
\end{proof}

%%% Local Variables: 
%%% mode: latex
%%% Te\Xx-master: "deform"
%%% Te\Xx-master: "deform"
%%% TeX-master: "deform"
%%% End: 

\vskip1cm
\subsection{A Factorization of the Pullback}\label{sec:pullback}
Let $L$ be a closed subgroup in a connected semisimple Lie group $G$ with
finite center and associated symmetric space $\Xx$, and let $\rho:L\to G'$ 
be a continuous homomorphism into a topological group $G'$.  
Combining the diagram in (\ref{eq:transfer-L})
with pullbacks in ordinary and bounded cohomology, 
we obtain the following commutative diagram:
\begin{equation}\label{eq:little-infty}
\xymatrix{
 \hcb^\bullet(G',\RR)
       \ar[rr]^{{\tt c}^\bullet_{G'}}\ar[d]_{\rho^\bullet_\mathrm{b}}
&
&\hc^\bullet(G',\RR)\ar[d]^{\rho^\bullet}
& \\
 \hb^\bullet(L,\RR)
       \ar[rr]^{{\tt c}^\bullet_L}\ar[rrd]_{\delta_{\infty,L}^\bullet}
&
&\h^\bullet(L,\RR)\ar[r]^\cong
&\h^\bullet\big(\Omega^\bullet(\Xx)^L\big)\\
&
&\h^\bullet\big(\Omega^\bullet_\infty(\Xx)^L\big)\ar[ur]_{i^\bullet_{\infty,L}}
& 
}
\end{equation}
from which one immediately reads:

\begin{cor}\label{cor:thm1-L} Let $G'$ be a topological group, $L\leq G$ 
any closed subgroup in a semisimple Lie group $G$ with finite center
and associated symmetric space $\Xx$ and $\rho:L\to G'$ a continuous 
homomorphism.  If $\alpha\in\hc^n(G',\RR)$ is represented by a continuous 
bounded class, then $\rho^{(n)}(\alpha)\in\h^n(L,\RR)$ is representable
by a $L$-invariant smooth closed differential $n$-form on $\Xx$
which is bounded.
\end{cor}

Analogously, if in addition $L=\Gamma<G$ is a lattice, then combining 
the top part of the diagram in (\ref{eq:bounded-rep})
with pullbacks we obtain

\begin{equation}\label{eq:little-p}
\xymatrix{
 \hcb^\bullet(G',\RR)
       \ar[rr]^{{\tt c}^\bullet_{G'}}\ar[d]_{\rho^\bullet_\mathrm{b}}
&
&\hc^\bullet(G',\RR)\ar[d]^{\rho^\bullet}
& \\
 \hb^\bullet(\Gamma,\RR)
       \ar[rr]^{{\tt c}^\bullet_\Gamma}\ar[rrd]_{\delta_{p,\Gamma}^\bullet}
&
&\h^\bullet(\Gamma,\RR)\ar[r]^\cong
&\h^\bullet\big(\Omega^\bullet(\Xx)^\Gamma\big)\\
&
&\h^\bullet\big(\Omega^\bullet_p(\Xx)^\Gamma\big)\ar[ur]_{i^\bullet_{p,\Gamma}}
& 
}
\end{equation}
In this section we shall mainly draw consequences from this, 
in especially relevant circumstances.  For example, 
if $G'$ also is a connected, semisimple Lie group with finite center, 
then in degree two the comparison map 
\begin{equation*}
{\tt c}^{(2)}_{G'}:\hcb^2(G',\RR)\to\hc^2(G',\RR)
\end{equation*}
is an isomorphism \cite{Burger_Monod_GAFA}, and we may then compose
$\big({\tt c}^{(2)}_{G'}\big)^{-1}$ with $\rho_\mathrm{b}^{(2)}$ and
$\delta_{p,\Gamma}^{(2)}$ to get a map, which we denote
by $\rho_p^{(2)}$, for which the following holds:

\begin{cor} \label{cor:factor-lp} 
If $G,G'$ are connected semisimple Lie groups with finite center,
$\Xx$ is the symmetric space associated to $G$ and $\Gamma<G$ is a lattice,
then the pullback via the homomorphism $\rho:\Gamma\to G'$ 
in ordinary cohomology and in degree two factors via
$\lp$-cohomology
\begin{equation*}
\xymatrix{
 \hc^2(G',\RR)\ar[d]^{\rho^{(2)}}\ar@/_3pc/[dd]_{\rho^{(2)}_p}
& \\
 \h^2(\Gamma,\RR)\ar[r]^\cong
&\h^2\big(\Omega^\bullet(\Xx)^\Gamma\big)\\
 \h^2\big(\Omega^\bullet_p(\Xx)^\Gamma\big)\ar[ur]_{i^{(2)}_{p,\Gamma}}
& 
}
\end{equation*}
\end{cor}

\begin{rem} (1) This is true for all closed subgroups $L<G$ in the case 
$p=\infty$. 

(2) Notice that, so far, we have not used the commutativity of the lower
part of the diagram in (\ref{eq:bounded-rep}).  This will be done
in the following section, to identify a numerical invariant
associated to a representation.
\end{rem}

%%% Local Variables: 
%%% mode: latex
%%% TeX-master: "deform"
%%% TeX-master: "deform"
%%% End: 

\vskip1cm
\subsection{The Hermitian Case}\label{subsec:hermitian}
Let $\Xx'$ be the symmetric space associated to $G'$, assume that
both $\Xx$ and $\Xx'$ are Hermitian symmetric, 
and that moreover $\Xx$ is irreducible.
Let $\kappa_\Xx\in\hc^2(G,\RR)$ and
$\kappa_\Xx^\mathrm{b}\in\hcb^2(G,\RR)$ be the cohomology classes
which corresponds to the \kahler form 
$\omega_\Xx$ on $\Xx$ via the isomorphisms
\begin{equation*}
\begin{aligned}
\Omega^2(\Xx)^G&\cong\hc^2(G,\RR)\cong\hcb^2(G,\RR)\\
\omega_\Xx\quad&\leftrightarrow\quad\kappa_\Xx\quad\,\,\leftrightarrow\quad\,\,\,\kappa_\Xx^\mathrm{b}\,.
\end{aligned}
\end{equation*}
The fact that $\hcb^2(G,\RR)\cong\RR\cdot\kappa_\Xx^\mathrm{b}$ 
and the definition of 
\begin{equation*}
\mathrm{T}_\mathrm{b}^{(2)}(\rho): \hcb^2(G',\RR)\to\hcb^2(G,\RR)
\end{equation*}
lead to the definition of the {\it bounded Toledo invariant} 
$\mathrm{t}_\mathrm{b}(\rho)$ by 
\begin{equation}\label{eq:bti}
\mathrm{T}_\mathrm{b}^{(2)}(\rho)(\kappa_{\Xx'}^\mathrm{b})=
\mathrm{t}_\mathrm{b}(\rho)\kappa_\Xx^\mathrm{b}\,.
\end{equation}
Then we have a Milnor--Wood type inequality:

\begin{lemma}\label{lem:milnor-wood} With the above notations,
\begin{equation*}
|\mathrm{t}_\mathrm{b}(\rho)|\leq\frac{\operatorname{rk}\Xx'}{\operatorname{rk}\Xx}\,.
\end{equation*}
\end{lemma}
\begin{proof} This follows from the fact that for any Hermitian symmetric space $\Yy$
with the appropriate normalization of the metric (see \S~\ref{sec:prelim})
\begin{equation*}
\|\kappa_\Yy^\mathrm{b}\|=\pi\operatorname{rk}\Yy\,,
\end{equation*}
and the fact that $\mathrm{T}_\mathrm{b}^\bullet(\rho)$ is norm decreasing 
in bounded cohomology.
\end{proof}

%{\bf DO NOT ERASE FOR THE MOMENT MY BEAUTIFUL DIAGRAM!!!}
%\begin{equation*}
%\xymatrix{
% \hcb^\bullet(G',\RR)
%       \ar[rr]^{{\tt c}^\bullet_{G'}}\ar[d]_{\rho^\bullet_\mathrm{b}}
%&
%&\hc^\bullet(G',\RR)\ar[d]^{\rho^\bullet}\ar@/_15pc/[dd]^{\rho^\bullet_p}
%& \\
% \hb^\bullet(\Gamma,\RR)
%       \ar[rrd]_{\delta_p^\bullet}
%&
%&\h^\bullet(\Gamma,\RR)\ar[r]^\cong
%&\h^\bullet(\Omega^\bullet(\Hh_\CC^n)^\Gamma)\\
%&
%&\h^\bullet(\Omega^\bullet_p(\Hh_\CC^n)^\Gamma)\ar[ur]_{i^\bullet_p}
%& 
%}
%\end{equation*}

The bounded Toledo invariant can now be nicely interpreted using
the lower part of (\ref{eq:bounded-rep}) in the case $p=2$.
In fact, the space $\Xx$ being Hermitian symmetric,
the $\ltwo$-cohomology spaces $\h^\bullet\big(\Omega_2^\bullet(\Xx)^\Gamma\big)$
are reduced and finite dimensional, \cite{Borel_64}[\S~3]. 
The following observation will be essential:

\begin{lemma}\label{lem:orth} Let $\Xx$ be a Hermitian symmetric space and $\Gamma$ 
a lattice in the isometry group $G:=\Iso(\Xx)^\circ$. 
Then the map 
\begin{equation*}
j_2^\bullet:\h^\bullet\big(\Omega^\bullet_2(\Xx)^\Gamma\big)\to
               \h^\bullet\big(\Omega^\bullet(\Xx)^G\big)=\Omega^\bullet(\Xx)^G
\end{equation*}
is the orthogonal projection.
\end{lemma}

\begin{proof} Denoting by $\<\,\cdot\,,\,\cdot\,\>_x$ the scalar product
on $\Lambda^\bullet(T_x\Xx)^\ast$, the scalar product
of two forms $\alpha,\beta\in\Omega_2^\bullet(\Xx)^\Gamma$
is given by
\begin{equation}\label{eq:scalar}
\<\alpha,\beta\>:=\int_{\Gamma\backslash\Xx}\<\alpha_x,\beta_x\>_xdv(\dot x)\,,
\end{equation}
where $dv$ is the volume measure on $\Gamma\backslash\Xx$;
fixing $x_0\in\Xx$, and letting $\mu$ be the $G$-invariant
probability measure on $\Gamma\backslash G$,
(\ref{eq:scalar}) can be written as
\begin{equation}\label{eq:scalar2}
 \<\alpha,\beta\>
=\vol(\Gamma\backslash G)
 \int_{\Gamma\backslash G}\<\alpha_{hx_0},\beta_{hx_0}\>_{hx_0}d\mu(\dot h)\,.
\end{equation}
Since we have identified $\h^\bullet\big(\Omega_2(\Xx)^\Gamma\big)$ with
the space of harmonic forms which are $\mathrm{L}^2$ (modulo $\Gamma$),
it suffices to show that 
\begin{equation*}
\big\<j^\bullet_2(\alpha),\beta\big\>=\big\<\alpha,j^\bullet_2(\beta)\big\>\,.
\end{equation*}
To this end we compute
\begin{equation*}
 \big\<(L_g^\ast\alpha)_x,\beta_x\big\>_x
=\<\alpha_{gx}\circ\Lambda^\bullet d_xL_g,\beta_x\>_x
=\big\<\alpha_{gx},\beta_x\circ(\Lambda^\bullet d_xL_g)^{-1}\big\>_{gx}
\end{equation*}
and hence, using (\ref{eq:scalar2}),
\begin{equation*}
\begin{aligned}
 &\<j^\bullet_2(\alpha),\beta\>\\
=&\vol(\Gamma\backslash G)\int_{\Gamma\backslash G}
    \left(\int_{\Gamma\backslash G}
         \<\alpha_{ghx_0},\beta_{hx_0}\circ(\Lambda^\bullet d_{hx_0}L_g)^{-1}\>_{ghx_0}d\mu(\dot g)\right)d\mu(\dot h)\\
=&\vol(\Gamma\backslash G)\int_{\Gamma\backslash G}
    \left(\int_{\Gamma\backslash G}
         \<\alpha_{gx_0},\beta_{hx_0}\circ(\Lambda^\bullet d_{hx_0}L_{gh^{-1}})^{-1}\>_{gx_0}d\mu(\dot g)\right)d\mu(\dot h)\\
=&\vol(\Gamma\backslash G)\int_{\Gamma\backslash G}
     \bigg\<\alpha_{gx_0},\int_{\Gamma\backslash G}\beta_{hx_0}\circ
        (\Lambda^\bullet d_{hx_0}L_{gh^{-1}})^{-1}d\mu(\dot h)\bigg\>_{gx_0}d\mu(\dot g)\,.
\end{aligned}
\end{equation*}
But $(\Lambda^\bullet d_{hx_0}L_{gh^{-1}})^{-1}=
     \Lambda^\bullet d_{gx_0}L_{hg^{-1}}$, so
\begin{equation*}
\begin{aligned}
\int_{\Gamma\backslash G}\beta_{hx_0}\circ
        (\Lambda^\bullet d_{hx_0}L_{gh^{-1}})^{-1}d\mu(\dot h)
&=\int_{\Gamma\backslash G}
        \beta_{hx_0}\circ \Lambda^\bullet d_{gx_0}L_{hg^{-1}}d\mu(\dot h)\\
&=\int_{\Gamma\backslash G}\beta_{hgx_0}\circ\Lambda^\bullet d_{gx_0}L_h d\mu(\dot h)
\end{aligned}
\end{equation*}
and hence, using (\ref{eq:scalar2}) and (\ref{eq:scalar},
\begin{equation*}
 \big\<j^\bullet _2(\alpha),\beta\big\>
=\int_{\Gamma\backslash\Xx}
    \left\<\alpha_x,\int_{\Gamma\backslash G}\beta_{hx}\circ
       \Lambda^\bullet d_x L_h d\mu(\dot h)\right\>_x dv(\dot x)
= \big\<\alpha,j^\bullet_2(\beta)\big\>
\end{equation*}
which shows that $j_2$ is self-adjoint.  Being clearly a projection, 
this proves the lemma.
\end{proof}

Applying the lemma in degree two, we have that 
for $\alpha\in\h^2\big(\Omega_2^\bullet(\Xx)^\Gamma\big)$,
\begin{equation*}
j^{(2)}_2(\alpha)=\frac{\<\alpha,\omega_\Xx\>}{\<\omega_\Xx,\omega_\Xx\>}\omega_\Xx\,,
\end{equation*}
where $\omega_\Xx$ is, as usual, the \kahler form on $\Xx$,
which is a generator of $\h^2\big(\Omega_2^\bullet(\Xx)^G\big)$
since $\Xx$ is assumed to be irreducible.
Define now
\begin{equation}\label{eq:ipi}
i_\rho:=\frac{\<\rho^{(2)}_2(\kappa_{\Xx'}),\omega_\Xx\>}{\<\omega_\Xx,\omega_\Xx\>}\,.
\end{equation}
It finally follows from the commutativity of the diagram 
\begin{equation}\label{eq:big-one}
\xymatrix{
 \hcb^\bullet(G',\RR)
       \ar[rr]^{{\tt c}^\bullet_{G'}}\ar[d]_{\rho^\bullet_\mathrm{b}}\ar@/_3pc/[ddd]_{\mathrm{T}^\bullet_\mathrm{b}(\rho)}
&
&\hc^\bullet(G',\RR)\ar[d]^{\rho^\bullet}
& \\
 \hb^\bullet(\Gamma,\RR)\ar[dd]_{\mathrm{T}_\mathrm{b}^\bullet}
       \ar[rr]^{{\tt c}^\bullet_\Gamma}\ar[rrd]_{\delta_{p,\Gamma}^\bullet}
&
&\h^\bullet(\Gamma,\RR)\ar[r]^\cong
&\h^\bullet\big(\Omega^\bullet(\Xx)^\Gamma\big)\\
&
&\h^\bullet\big(\Omega^\bullet_p(\Xx)^\Gamma\big)
   \ar[ur]_{i^\bullet_{p,\Gamma}}\ar[dr]^{j_p^\bullet}
& \\
 \hcb^\bullet(G,\RR)\ar[rr]^{{\tt c}^\bullet_G}
&
&\hc^\bullet(G,\RR)\ar[r]_\cong
&\Omega^\bullet(\Xx)^G\,.
}
\end{equation}
in the special case of $p=2$ and degree 2 and from Corollary~\ref{cor:factor-lp} that:

\begin{lemma}\label{lem:eq} $i_\rho=\mathrm{t}_\mathrm{b}(\rho)$.
\end{lemma}

For our immediate applications, we draw the following conclusions 
when $\Xx=\Hh_\CC^p$:

\begin{prop} Let $G$ be a connected simple Lie group with finite center
and associated symmetric space $\Hh_\CC^p$ and 
let $\Gamma<G$ be a lattice.
Let $\Xx'$ be a Hermitian symmetric space, $G':=\Iso(\Xx')^\circ$ and
$\rho:\Gamma\to G'$ a representation.  Then:
\begin{enumerate}
\item $|i_\rho|\leq\operatorname{rk}\Xx'$;
\medskip
\item if either $\Gamma\backslash\Hh_\CC^p$ is compact or $n\geq2$,
then $i_\rho$ is a characteristic number.
\end{enumerate}
\end{prop}

\begin{proof} (1) follows from Lemma~\ref{lem:milnor-wood} and Lemma~\ref{lem:eq}.

\medskip
\noindent
(2) If $\Gamma\backslash\Hh_\CC^p$ is compact we have evidently that
\begin{equation*}
\h^\bullet\big(\Omega_2^\bullet(\Hh_\CC^p)^\Gamma\big)=
\h^\bullet\big(\Omega^\bullet(\Hh_\CC^p)^\Gamma\big)\,,
\end{equation*}
while if  $p\geq2$, we have that  $\h^2\big(\Omega_2^\bullet(\Hh_\CC^p)^\Gamma\big)$
injects into $\h^2\big(\Omega^\bullet(\Hh_\CC^p)^\Gamma\big)$ \cite{Zucker}.
In any case, by Corollary~\ref{cor:factor-lp},  $i_\rho$ depends only 
upon $\rho^{(2)}(\kappa_{\Xx'})\in\h^2\big(\Omega^\bullet(\Hh_\CC^p)^\Gamma\big)$, 
which is well known to be a characteristic class (see \S~\ref{sec:intro}).  
Hence $i_\rho$ is a characteristic number.
\end{proof}

Actually, the injectivity of the map 
$\h^2\big(\Omega_2^\bullet(\Hh_\CC^p)^\Gamma\big)
\to
\h^2\big(\Omega^\bullet(\Hh_\CC^p)^\Gamma\big)$
is proven in \cite{Zucker} in a more general framework but only 
for arithmetic lattices.
In our specific case however the proof applies 
without the arithmeticity requirement \cite{Saper_email}.
%%% Local Variables: 
%%% mode: latex
%%% TeX-master: "deform"
%%% End: 

\vskip1cm
\subsection{The Formula, Once Again}\label{sec:formula}
The very explicit form of the factorization of the comparison map
between bounded and ordinary cohomology, 
together with the implementation of the pullback
by boundary maps in \cite{Burger_Iozzi_app}
allows one to give explicit representatives of the class
$\rho^{(2)}(\kappa_q)$ at least when $\Xx'=\Hh_\CC^q$.  

\bigskip

To this end we assume that 
the homomorphism $\rho:L\to\Iso(\Hh_\CC^q)^\circ=:G'$ is nonelementary
so that there exists a $L$-equivariant measurable map   
$\varphi:\Hh_\CC^p\to\Hh_\CC^q$ (see for example \cite{Burger_Mozes}).  Let
\begin{equation*}
c_q:(\partial\Hh_\CC^q)^3\to[-1,1]
\end{equation*}
be the Cartan cocycle which is a representative of the \kahler class
$\kappa_q^\mathrm{b}\in\hcb^2(G',\RR)$, \cite{Burger_Iozzi_supq}.
For $\xi\in\partial\Hh_\CC^\ell$, 
let $e^\xi$ denote the exponential of the Busemann function defined in
(\ref{eq:exi}).  Then

\begin{prop}\label{prop:bounded-rep}  Let $G$ be a connected simple Lie group 
with finite center and associated symmetric space $\Hh_\CC^p$ and let 
$L\leq G$ be any closed subgroup.  Then the differential $2$-form
\begin{equation}\label{eq:explicit+} 
\int_{(\partial\Hh_\CC^p)^3} 
      c_q\big(\varphi(\xi_0),\varphi(\xi_1),\varphi(\xi_2)\big)
   e^{\xi_0}\wedge de^{\xi_1}\wedge de^{\xi_2} 
      d\nu_0^3(\xi_0,\xi_1,\xi_2)
\end{equation}
is a smooth $L$-invariant bounded closed $2$-form representing 
$\rho^{(2)}(\kappa_q)\in\h^2(L,\RR)\cong\h^2\big(\Omega^\bullet(\Hh_\CC^p)^L\big)$.
\end{prop}

\begin{proof} Since $c_q\in\Bb^\infty\big((\partial\Hh_\CC^q)^3\big)^{G'}$
represents $\kappa_q^\mathrm{b}$, 
by (\ref{eq:old-rep}) the cocycle in $\linfty\big((\partial\Hh_\CC^p)^3\big)^L$
\begin{equation*}
(\xi_0,\dots,\xi_n)\mapsto c\big(\varphi(\xi_0),\dots,\varphi(\xi_n)\big)
\end{equation*}
represents canonically 
$\rho^{(2)}_\mathrm{b}(\kappa_q^\mathrm{b})\in\hb^2(L,\RR)$.
By Lemma~\ref{lem:explicit}, (\ref{eq:explicit+}) is a smooth
differential $2$-form in $\Omega_\infty^2(\Xx)$ which is
$L$-invariant and, by Proposition~\ref{prop:transfer-L},
it represents $\rho^{(2)}(\kappa_q)\in\h^2\big(\Omega^\bullet(\Hh_\CC^p)^L\big)$.
\end{proof}

Let us assume now that $L=\Gamma<\SUp$ is a lattice and move to the main formula, 
which will be an implementation
of \S~\ref{subsec:formula} in our concrete situation.
Let  $\Cc_p$ be the set of all chains in $\partial\Hh_\CC^p$
and, for any $k\geq1$, let 
\begin{equation*}
\Cc_p^{\{k\}}:=\big\{(C,\xi_1,\dots,\xi_k):\,C\in\Cc_p,(\xi_1,\dots,\xi_k)\in C^k\big\}
\end{equation*}  
be the space of configurations of $k$-tuples of points on a chain.  
Both $\Cc_p$ and $\Cc_p^{\{1\}}$ are homogeneous spaces of $\SUp$.
In fact, the stabilizer $H$ in $G$ of a fixed chain $C_0\in\Cc_p$
is also the stabilizer of a plane of signature
$(1,1)$ in $\SUp$ and hence isomorphic to 
$\mathrm{S}\big(\mathrm{U(1,1)}\times\mathrm{U(p-1)}\big)$.
Then $\SUp$ acts transitively on $\Cc_p$ (for example because it acts
transitively on pairs of points in $\partial\Hh_\CC^p$ and any
two points in $\partial\Hh_\CC^p$ determine uniquely a chain)
and $H$ acts transitively on $C_0$, so that, if $Q=P\cap H$,
where $P$ is the stabilizer in $\SUp$ of a fixed basepoint $\xi_0\in C_0$,
there are $\SUp$-equivariant 
(hence measure class preserving) diffeomorphisms
\begin{equation*}
\begin{aligned}
\SUp/H&\to\,\,\Cc_p\\
gH\,\,&\mapsto gC_0
\end{aligned}
\end{equation*}
and
\begin{equation*}
\begin{aligned}
\SUp/Q&\to\quad\,\Cc_p^{\{1\}}\\
gQ\,\,\,&\mapsto(gC_0,g\xi_0)\,.
\end{aligned}
\end{equation*}
Moreover, the projection $\pi:\Cc_p^{\{1\}}\to\Cc_p$
which associates to a point $(C,\xi)\in\Cc_p^{\{1\}}$ 
the chain $C\in\Cc_p$ is a $\SUp$-equivariant fibration,
the space $\Cc_p^{\{k\}}$ appears then naturally as 
$k$-fold fibered product of $\Cc_p^{\{1\}}$ with respect to $\pi$,
and for every $k\geq1$, the map
\begin{equation}\label{eq:chains-fibered}
\begin{aligned}
(\SUp/Q)_f^k\quad\,\,\,&\to\qquad\quad\Cc_p^{\{k\}}\\
(x_1Q,\dots,x_kQ)&\mapsto(gC_0,x_1\xi_0,\dots,x_k\xi_0)
\end{aligned}
\end{equation}
where $x_iH=gH$, $1\leq i\leq k$, is a $\SUp$-equivariant diffeomorphism
which preserves the $\SUp$-invariant Lebesgue measure class.
Using Fubini's theorem, one has that for almost every $C\in\Cc_p$
the restriction 
\begin{equation*}
\varphi_C:C\to\partial\Hh_\CC^q
\end{equation*}
of $\varphi$ to $C$ 
is measurable and for every $\gamma\in \Gamma$ and almost every $\xi\in C$
\begin{equation*}
\varphi_{\gamma C}(\gamma \xi)=\rho(\gamma)\varphi_C(\xi)\,.
\end{equation*}
This allows us to define
\begin{equation*}
\begin{aligned}
\varphi^{\{3\}}:\qquad\Cc_p^{\{3\}}\,\,\,\,\,\,\,\,&\to\qquad\qquad(\partial\Hh_\CC^q)^3\\
(C,\xi_1,\xi_2,\xi_3)&\mapsto\big(\varphi_C(\xi_1),\varphi_C(\xi_2),\varphi_C(\xi_3)\big)\,.
\end{aligned}
\end{equation*}
Then
\begin{thm}\label{thm:formula} 
Let $i_\rho$ be the invariant defined in (\ref{eq:ipi}).
Then for almost every chain $C\in\Cc_p$ and almost every $(\xi_1,\xi_2,\xi_3)\in C^3$,
\begin{equation*}
\int_{L\backslash\SUp}c_q\big(\varphi^{\{3\}}(gC,g\xi_1,g\xi_2,g\xi_3)\big)d\mu(g)=
  i_\rho c_p(\xi_1,\xi_2,\xi_3)\,,
\end{equation*}
where $c_\ell$ is the Cartan invariant
and $\mu$ is the $\SUp$-invariant probability measure on 
$\Gamma\backslash\SUp$.
\end{thm}

\begin{cor}\label{cor:formula} Assume that $i_\rho=1$.  Then for almost every 
$C\in\Cc_p$ and almost every $(\xi_1,\xi_2,\xi_3)\in C^3$
\begin{equation*}
c_q\big(\varphi_C(\xi_1),\varphi_C(\xi_2),\varphi_C(\xi_3)\big)=c_p(\xi_1,\xi_2,\xi_3)\,.
\end{equation*}
\end{cor}

\begin{proof}[Proof of Theorem~\ref{thm:formula}]
Let $H,P,Q<\SUp$ such as in the above discussion.  Since $P$
is the stabilizer of a basepoint $\xi_0\in\partial\Hh_\CC^p$,
it is a minimal parabolic subgroup and hence the closed subgroup
$Q$ is amenable.  Moreover, $H$ acts ergodically on $H/Q\times H/Q$
since in $H/Q\times H/Q$ there is an open $H$-orbit of full measure.
We can hence apply Proposition~\ref{prop:formula-homog}
with $G=\SUp$, $G'=\mathrm{PU}(q,1)$ and $\kappa'=\kappa_q^\mathrm{b}$.
Moreover, by (\ref{eq:bti}) and the definition of the bounded Toledo map, 
let $\kappa=\mathrm{T}_\mathrm{b}^{(2)}\big(\rho_\mathrm{b}^{(2)}(\kappa_q^\mathrm{b})\big)=\mathrm{t}_\mathrm{b}(\rho)\kappa_p^\mathrm{b}$,
which in turns, by Lemma~\ref{lem:eq} implies that 
$\kappa=i_\rho\kappa_p^\mathrm{b}$.  Set $G/P=\partial\Hh_\CC^p$,
$c'=i_\rho c_p\in\Bb^\infty\big((\partial\Hh_\CC^p)^3\big)^{\SUp}$,
$X=\partial\Hh_\CC^p$ and 
$c'=c_q\in\Bb^\infty\big((\partial\Hh_\CC^p)^3\big)^{\mathrm{PU}(q,1)}$.
Then the conclusion of the theorem is immediate if we observe that 
the identification in (\ref{eq:chains-fibered}) transforms the
map $\varphi^3_f$ defined in (\ref{eq:varphifn}) into 
the map $\varphi^{\{3\}}$ defined above.
\end{proof}

\begin{rem}\label{rem:fibered}  
It is now clear what is the essential use of the fibered
product: the triples of points that lie on a chain form a set of
measure zero in $\big(\partial\Hh_\CC^p\big)^3$, and hence we
would not have gained any information on these configuration
of points by the direct use of the more familiar formula as in 
Corollary~\ref{cor:usual-formula}.
\end{rem}

%\begin{exo} In our basic example with $G,H,Q,P$ as in \S~\ref{subsec:gpullback}
%and $c_p:(\partial\Hh_\CC^p)^3\to[-1,1]$ the Cartan cocycle, 
%we have that $c_p$ represents the bounded \kahler class 
%$\kappa_p^\mathrm{b}\in\hcb^2(G,\RR)$.
%Thus the cocycle $\tilde c_p:\Cc_p^{(3)}\to[-1,1]$
%defined by
%\begin{equation*}
%\tilde c_p(C,x_1,x_2,x_3)=
%\begin{cases}
%\hphantom{-}1 &\hbox{ if }x_1,x_2,x_3\hbox{ are positively oriented on }C\\
%-1            &\hbox{ if }x_1,x_2,x_3\hbox{ are negatively oriented on }C\\
%\hphantom{-}0 &\hbox{ if }x_1,x_2,x_3\hbox{ are not all distinct}\,,
%\end{cases}
%\end{equation*}
%is a representative of 
%$\alpha_G(\kappa_p^\mathrm{b})\in\hcb^2(G,\linfty(\Cc_n))$.
%\end{exo}

%%% Local Variables: 
%%% mode: latex
%%% Te\Xx-master: "deform"
%%% TeX-master: "deform"
%%% End: 

\vskip1cm
\section{The Measurable Cartan Theorem}\label{sec:cartan}
The goal of this section if to prove the following measurable
version of a theorem of E.~Cartan.

\begin{thm}\label{thm:cartan} Let $p\geq2$ and let
$\varphi:\partial\Hh_\CC^p\to\partial\Hh_\CC^q$
be a measurable map such that:
\begin{enumerate}
\item[(i)] for almost every chain $C$ and almost every triple $(\xi,\eta,\zeta)$
of distinct points on $C$, the triple 
$\varphi(\xi), \varphi(\eta), \varphi(\zeta)$ consists also of distinct
points which lie on a chain and have the same orientation as $(\xi,\eta,\zeta)$;
\item[(ii)] for almost every triple of points $\xi,\eta,\zeta$ not on a chain, 
$\varphi(\xi)$, $\varphi(\eta)$, $\varphi(\zeta)$ are also not on a chain.
\end{enumerate}
Then there is an isometric holomorphic embedding $F:\Hh_\CC^p\to\Hh_\CC^q$
such that $\partial F$ coincides with $\varphi$ almost everywhere.
\end{thm}

The proof goes as follows. We first show by induction that the
statement of the theorem for a fixed $p$ follows from the analogous
statement in one lower dimension, provided $p\geq3$;
this leaves us to show the statement for $p=2$.
The next step is to show that if $p=2$ any map
$\varphi:\partial\Hh_\CC^2\to\partial\Hh_\CC^q$ satisfying
the hypotheses of Theorem~\ref{thm:cartan} takes values
also in $\partial\Hh_\CC^2$; this will be achieved 
by an appropriate convex hull argument.
The last step is hence to show the assertion for $p=q=2$,
for which we need a careful modification of Cartan's argument
for pointwise defined maps.  

\subsection{Reduction to the case $p=2$}
If $k\leq p$, let $\Pp_k$ denote the set of $k$-planes (see the end
of \S~\ref{sec:prelim}) and, if $x\in\Hh_\CC^p$,
$\Pp_k(x)$ the subset of $\Pp_k$ of $k$-planes through $x$.
%For $W\in\Pp_k$,
%the boundary $\partial W\subset\Hh_\CC^p$ is called a {\it $k$-chain}
%and the set of $k$-chains is denoted by $\Cc_k$.  

We now let $p\geq3$, we assume that the theorem holds for $p-1$
and we want to show that then it holds for $p$.  
Let us start by observing that a simple verification using Fubini's theorem 
applied to the configuration spaces
\begin{equation*}
\big\{(X,C):\,X\in\Pp_{p-1},C\in\Cc_p, C\subset\partial X\big\}
\end{equation*}
and
\begin{equation*}
\big\{(X,\xi_1,\xi_2,\xi_3):\,X\in\Pp_{p-1},\,\,\,\text{and }\xi_1,\xi_2,\xi_3\in\partial X\big\}
\end{equation*}
shows that, for almost every $X\in\Pp_{p-1}$, the restriction
$\varphi|_{\partial X}$ of $\varphi$ to $\partial X$
is measurable and satisfies the hypotheses of Theorem~\ref{thm:cartan}.

Applying the induction hypothesis we get for almost every $X\in\Pp_{p-1}$ 
an isometric holomorphic embedding
\begin{equation*}
F_X:X\to \Hh_\CC^q
\end{equation*}
such that $\partial F_X=\varphi|_{\partial X}$ almost everywhere.
Thus the set
\begin{equation*}
\begin{aligned}
\big\{
(x,X):\,X\in\Pp_{p-1}(x) \hbox{ and  there is }F_X:X\to\Hh_\CC^q
      \hbox{ as above}&\\
 \hbox{ with }\partial F_X=\varphi|_{\partial X}\hbox{ almost everywhere}&
\big\}
\end{aligned}
\end{equation*}
is of full measure in the configuration space
\begin{equation*}
\big\{(x,X):X\in\Pp_{p-1}(x)\big\}\,,
\end{equation*}
and we may define for almost every $x\in\Hh_\CC^p$ and 
almost every $X\in\Pp_{p-1}(x)$ the function
\begin{equation*}
f(x,X):=F_X(x)\,.
\end{equation*}

Using again Fubini's theorem, one checks that for almost every 
$X_1$, $X_2\in\Pp_{p-1}(x)$,
\begin{equation*}
\partial F_{X_1}|_{\partial X_1\cap \partial X_2}=
   \varphi|_{\partial X_1\cap\partial X_2}=
   \partial F_{X_2}|_{\partial X_1\cap\partial X_2}
\end{equation*}
and, since $p\geq3$,
\begin{equation*}
\partial X_1\cap\partial X_2=\partial(X_1\cap X_2)\neq\emptyset\,,
\end{equation*}
which implies that $f(x,X_1)=f(x,X_2)$.
Thus $f(x,X)$ is almost everywhere independent of $X\in\Pp_{p-1}(x)$
and gives rise to a well defined map $f:\Hh_\CC^p\to\Hh_\CC^q$ 
which by construction preserves the distances of almost every pair 
of points.  It is then not difficult to see that $f$ coincides
almost everywhere with an isometric embedding $\Hh_\CC^p\to\Hh_\CC^q$.
This, together with the fact that $\partial f=\varphi$
preserves the orientation on chains, implies that 
the embedding must be holomorphic.
\qed

\subsection{Reduction to the case $p=q=2$}
Denote by $(\partial\Hh_\CC^p)^{(k)}$ the subset of full measure in 
$(\partial\Hh_\CC^p)^k$ consisting of $k$-tuples of distinct points
in $\partial\Hh_\CC^p$. Recall that any two distinct chains are 
either disjoint or intersect in a point,
and hence every pair of distinct points $(\xi,\eta)\in(\partial\Hh_\CC^p)^{(2)}$
determines a unique chain $C(\xi,\eta)$.

\begin{lemma}\label{lem:ch-pres} 
Let $\varphi:\partial\Hh_\CC^p\to\partial\Hh_\CC^q$
be a measurable map satisfying the hypothesis (i) of Theorem~\ref{thm:cartan}
and let $c_\ell:(\partial\Hh_\CC^\ell)^3\to[-1,1]$ be the Cartan cocycle.
Then:
\begin{enumerate}
\item for almost every $(\xi_1,\xi_2)\in(\partial\Hh_\CC^p)^{(2)}$,
we have that $\varphi(\xi_1)\neq\varphi(\xi_2)$, and
\item for almost every $\xi_3\in C(\xi_1,\xi_2)$, 
  we have that 
\begin{equation*}
\varphi(\xi_3)\in C\big(\varphi(\xi_1),\varphi(\xi_2)\big)
\end{equation*}
  and 
\begin{equation*}
c_q\big(\varphi(\xi_1),\varphi(\xi_2),\varphi(\xi_3)\big)=c_p(\xi_1,\xi_2,\xi_3)\,.
\end{equation*}
\end{enumerate}
\end{lemma}

As a consequence we have: 

\begin{cor}\label{cor:big-fhi} 
Let, as above, $\varphi:\partial\Hh_\CC^p\to\partial\Hh_\CC^q$
be a measurable map satisfying the hypothesis (i) of Theorem~\ref{thm:cartan}
and let $c_\ell$ be the Cartan cocycle.
Then there is a measurable map 
\begin{equation}\label{eq:big-fhi}
\Phi:\Cc_p\to\Cc_q
\end{equation} 
such that 
\begin{equation}\label{eq:5.2}
\Phi\big(C(\xi_1,\xi_2)\big)=C\big(\varphi(\xi_1),\varphi(\xi_2)\big)
\end{equation}
for almost every $(\xi_1,\xi_2)\in(\partial\Hh_\CC^p)^2$.
\end{cor}

\begin{proof}[Proof of Lemma~\ref{lem:ch-pres}] 
Consider the measure class preserving bijection
\begin{equation*}
\begin{aligned}
(\partial\Hh_\CC^p)^{(2)}&
   \to\big\{(C,\xi_1,\xi_2):\,C\in\Cc_p,\xi_1,\xi_2\in C,\xi_1\neq \xi_2\big\}\\
(\xi_1,\xi_2)\,\,\,&\mapsto\qquad\qquad\big(C(\xi_1,\xi_2),\xi_1,\xi_2\big)\,.
\end{aligned}
\end{equation*}
Then Theorem~\ref{thm:cartan}(i) implies by Fubini that
for almost every $C\in\Cc_p$, for almost every $(\xi_1,\xi_2)\in C^{(2)}$ 
and for almost every $\xi_3\in C$ we have
\begin{equation}\label{eq:ch-pres}
c_q\big(\varphi(\xi_1),\varphi(\xi_2),\varphi(\xi_3)\big)=c_p(\xi_1,\xi_2,\xi_3)
\end{equation}
which, using the above bijection, is equivalent to the fact that
for almost every $(\xi_1,\xi_2)\in(\partial\Hh_\CC^p)^{(2)}$ and
for almost every $\xi_3\in C(\xi_1,\xi_2)$, (\ref{eq:ch-pres}) holds,
which shows that $\varphi(\xi_1)\neq\varphi(\xi_2)$ and 
that (2) holds.
\end{proof}

\begin{proof}[Proof of Corollary~\ref{cor:big-fhi}] It is clear that if $C\in\Cc_p$ is such that 
for almost every $(\xi_1,\xi_2)\in C^{(2)}$ and 
for almost every $\xi_3\in C$ (\ref{eq:ch-pres}) holds,
then in particular if $(\xi_1,\xi_2)\in C^{(2)}$ and $(\eta_1,\eta_2)\in C^{(2)}$
are such that 
\begin{equation*}
\begin{aligned}
c_q\big(\varphi(\xi_1),\varphi(\xi_2),\varphi(\xi_3)\big)&=c_p(\xi_1,\xi_2,\xi_3)\\
c_q\big(\varphi(\eta_1),\varphi(\eta_2),\varphi(\eta_3)\big)&=c_p(\eta_1,\eta_2,\eta_3)
\end{aligned}
\end{equation*}
then $\varphi(\xi_1)\neq\varphi(\xi_2)$, $\varphi(\eta_1)\neq\varphi(\eta_2)$
and  $C\big(\varphi(\xi_1),\varphi(\xi_2)\big)\cap C\big(\varphi(\eta_1),\varphi(\eta_2)\big)$
contains the essential image of $\varphi|_C$.
Since this cannot be reduced to a point, we have that 
$C\big(\varphi(\xi_1),\varphi(\xi_2)\big)=
C\big(\varphi(\eta_1),\varphi(\eta_2)\big)$,
which leads to the map
\begin{equation*}
\begin{aligned}
\Phi:\qquad\Cc_p\qquad&\to\qquad\Cc_q\\
\big(C(\xi_1,\xi_2)\big)&\mapsto C\big(\varphi(\xi_1),\varphi(\xi_2)\big)
\end{aligned}
\end{equation*}
which is then well defined and satisfies (\ref{eq:5.2}).
\end{proof}

Now choose a Borel map
\begin{equation*}
\begin{aligned}
(\partial\Hh_\CC^p)^{(2)}&\to\Mm^1(\partial\Hh_\CC^p)\\
(\xi_1,\xi_2)\,\,&\mapsto\,\,\,\,\,\mu_{(\xi_1,\xi_2)}\,,
\end{aligned}
\end{equation*}
such that $\mu_{(\xi_1,\xi_2)}$ is in the Lebesgue measure class of $C(\xi_1,\xi_2)$.  
Let $d\lambda$ be the ``round measure'' on $\partial\Hh_\CC^p$
and consider the map
\begin{equation*}
\begin{aligned}
M:\,\,\,\,(\partial\Hh_\CC^p)^{(3)}\,\,
      &\to\Mm^1(\partial\Hh_\CC^p\times\partial\Hh_\CC^p)\\
(\xi_1,\xi_2,\xi_3)
      &\mapsto\,\mu_{(\xi_1,\xi_2)}\otimes\mu_{(\xi_1,\xi_3)}\,.
\end{aligned}
\end{equation*}

\begin{lemma}\label{lem:dir-im}
The measure on $\partial\Hh_\CC^p\times\partial\Hh_\CC^p$
defined by 
\begin{equation*}
\int_{(\partial\Hh_\CC^p)^{(3)}}\big(\mu_{(\xi_1,\xi_2)}\otimes\mu_{(\xi_1,\xi_3)}(f)\big)d\lambda^3(\xi_1,\xi_2,\xi_3)
\end{equation*}
for $f\in\mathrm{C}\big(\partial\Hh_\CC^p\times\partial\Hh_\CC^p\big)$,
is equivalent to $\lambda^2$.
\end{lemma}
\begin{proof} This is obvious.
\end{proof}

\begin{lemma}\label{lem:belong} For almost every 
$(\xi_1,\xi_2,\xi_3)\in(\partial\Hh_\CC^p)^{(3)}$,
for almost every $(a,b)\in C(\xi_1,\xi_2)\times C(\xi_1,\xi_3)$
and for almost every $c\in C(a,b)$, we have that:
\begin{enumerate}
\item[--] $\varphi(a)\in C\big(\varphi(\xi_1),\varphi(\xi_2)\big)$;
\item[--] $\varphi(b)\in C\big(\varphi(\xi_1),\varphi(\xi_3)\big)$;
\item[--] $\varphi(c)\in C\big(\varphi(a),\varphi(b)\big)$.
\end{enumerate}
\end{lemma}
\begin{proof} This follows from repeated applications of
Lemma~\ref{lem:ch-pres}(2) and Lemma~\ref{lem:dir-im}.
\end{proof}

\begin{cor} Assume that in Theorem~\ref{thm:cartan} we have 
$p=2$.  Then the essential image of $\varphi$ is contained 
in a $2$-chain.
\end{cor}

\begin{proof} Fix $(\xi_1,\xi_2,\xi_3)$ not on a chain, for which
Lemma~\ref{lem:belong} holds.  
Let $E\subset\partial\Hh_\CC^2$ be the set of $c\in\partial\Hh_\CC^2$
such that there are $(a,b)\in C(\xi_1,\xi_2)\times C(\xi_1,\xi_3)$ with 
$c\in C(a,b)$ and Lemma~\ref{lem:belong} holds for $a,b,c$.
Then $E\subset\partial\Hh_\CC^2$ is of full measure.
Moreover, 
\begin{equation*}
\begin{aligned}
\varphi(c)&\hbox{ is in the $\CC$-linear span of }\varphi(a)\hbox{ and }\varphi(b),\\
\varphi(a)&\hbox{ is in the $\CC$-linear span of }\varphi(\xi_1)\hbox{ and }\varphi(\xi_2),\hbox{ and}\\
\varphi(b)&\hbox{ is in the $\CC$-linear span of }\varphi(\xi_1)\hbox{ and }\varphi(\xi_3)\,,
\end{aligned}
\end{equation*}
so that for all $c\in E$, $\varphi(c)$ is in the $2$-chain
determined by the $3$-dimensional space
$\varphi(\xi_1)\oplus\varphi(\xi_2)\oplus\varphi(\xi_3)$.
\end{proof}

\subsection{The Case $p=q=2$}
So we are reduced to show Theorem~\ref{thm:cartan} for $p=q=2$.
The essential point of the proof will be the following:

\begin{prop}\label{prop:cartan} Let $g:\CC\to\CC$ be a measurable
map such that for almost every circle $S\subset\CC$, there is
a circle $\Gamma(S)\subset\CC$ such that:
\begin{enumerate}
\item[(i)] for almost every $z\in S$, $g(z)\in\Gamma(S)$;
\item[(ii)] for almost every $z_1,z_2,z_3\in S$ distinct, 
$g(z_1),g(z_2), g(z_3)\in\Gamma(S)$ are distinct and in the same cyclic order;
\item[(iii)] for almost every $z\in\CC$ the set 
\begin{equation*}
\{(S_1,S_2):\,z\in S_i,\, i=1,2,\,\Gamma(S_1)=\Gamma(S_2)\}
\end{equation*}
is of measure zero.
\end{enumerate}
Then $g$ coincides almost everywhere with an affine map 
$z\mapsto\lambda z+c$, where $\lambda\in\CC^\times$ and $c\in\CC$.
\end{prop}

\medskip
\noindent
{\it Proof.}
By Fubini's theorem, the set $E$ of $z\in\CC$ such that (i), (ii) and (iii)
hold
for almost every circle through $z$ is of full measure in $\CC$;
fixing $z\in E$ and composing with an affine map, we may assume that $g(z)=z$.
Conjugating $g$ with an inversion $i:\widehat\CC\to\widehat\CC$
through a circle with center $z$, we get a map $f:\CC\to\CC$
which induces a map $F:\Dd\to\Dd$ on the set $\Dd$ of affine $\RR$-lines in $\CC$
satisfying the following properties:
\begin{enumerate}
\item[(ii)'] for almost every $d\in\Dd$, there is $F(d)\in\Dd$ such that
for almost every $A,B,C\in d$ distinct, $f(A), f(B), f(C)$ 
are distinct and lie on $F(d)$;
\item[(iii)'] for almost every point $w\in\CC$, the set
\begin{equation*}
\big\{(d_1,d_2)\in\Dd(w):\,F(d_1)=F(d_2)\big\}
\end{equation*}
is of measure zero.
Here $\Dd(w)$ is the set of lines through $w$.
\end{enumerate}
Notice that in (ii)' one could have stated also a condition corresponding
to the preservation of the cyclic ordering.  However this is never used
for the maps $f$ and $F$ themselves and comes into play only at the very end 
of the proof when we return to the map $g$, so that we chose to ignore it.
\begin{lemma}\label{lem:anharmonic} Let 
\begin{equation*}
\Kk=
\big\{
(d,A,B,C):\,d\in\Dd,\hbox{ and } A,B,C\hbox{ are distinct points on }d
\big\}\,.
\end{equation*}
Then for almost every $(d,A,B,C)\in\Kk$:
\begin{enumerate}
\item[{\it (i)}] we have
\begin{equation*}
\big(F(d),f(A),f(B),f(C)\big)\in\Kk\,,
\end{equation*} and
\item[{\it (ii)}] $f$ preserves the property for $4$-tuples of points
on a line to be in {\it anharmonic position}: namely,  
if $[\,\cdot\,,\,\cdot\,,\,\cdot\,,\,\cdot\,]$ denotes 
the crossratio, and if $D\in d$ and 
$D'\in F(d)$ are points such that 
$[A,B,C,D]=-1$ and $\big[f(A),f(B),f(C),D'\big]=-1$, then 
\begin{equation*}
f(D)=D'\,.
\end{equation*}
\end{enumerate}
\end{lemma}

\begin{proof}
Let 
\begin{equation*}
\Kk_0=
\big\{
(d,A);\,d\in\Dd,\,A\in d
\big\}
\end{equation*}
and let 
\begin{equation*}
\begin{aligned}
\Kk_1=
\big\{
(d',d,A,B,C,M):\,(d,A,B,C)\in\Kk,\,(d',A)\in\Kk_0&\\
\hbox{ and }d'\neq d,\,M\notin d'\cup d&\big\}\,.
\end{aligned}
\end{equation*}

First we observe that properties (ii)' and (iii)' together with a repeated use
of Fubini's theorem imply that there is a subset $E\subset\Kk_1$ of full measure 
such that for all $(d',d,A,B,C,M)\in E$, then
\begin{equation*}
\big(F(d'),F(d),f(A),f(B),f(C),f(M)\big)\in\Kk_1\,.
\end{equation*} 

Now to every $(d',d,A,B,C,M)\in\Kk_1$ associate the picture in 
Figure~\ref{fig:cartan}
%\begin{center}
 \begin{figure}
 \includegraphics[scale=0.497]{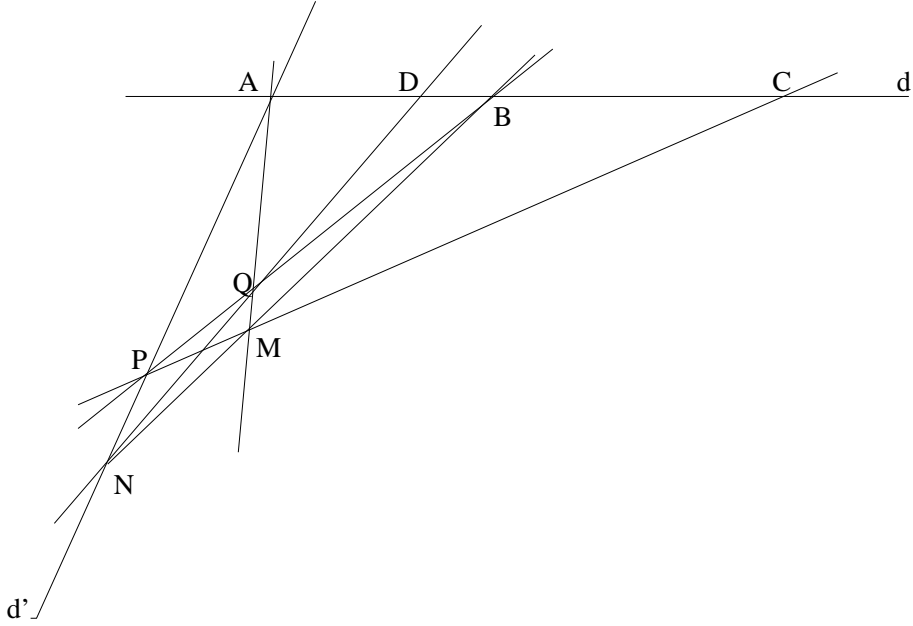}
\caption{{\it La m\'ethode du quadrilat\`ere complet}}
\label{fig:cartan}
\end{figure}
%\end{center}
and the following seven maps $m_i:\Kk_1\to\Kk$ given by
\begin{equation*}
\begin{aligned}
m_1:(d',d,A,B,C,M)&\mapsto\big(\<C,M\>,C,M,P\big)\\
m_2:(d',d,A,B,C,M)&\mapsto\big(\<P,B\>,P,Q,B\big)\\
m_3:(d',d,A,B,C,M)&\mapsto\big(\<A,M\>,A,Q,M\big)\\
m_4:(d',d,A,B,C,M)&\mapsto\big(\<B,M\>,B,M,N\big)\\
m_5:(d',d,A,B,C,M)&\mapsto(d',A,P,N)\\
m_6:(d',d,A,B,C,M)&\mapsto\big(\<N,Q\>,N,Q,D\big)\\
m_7:(d',d,A,B,C,M)&\mapsto(d,A,B,D)\,.
\end{aligned}
\end{equation*}
One verifies easily that 
for every set $T\subset\Kk$ of full measure,
$m_i^{-1}(T)\subset\Kk_1$ is of full measure for $1\leq i\leq7$.
So let $T\subset\Kk$ be the subset consisting of $(d,A,B,C)$
such that $\big(f(A),f(B),f(C)\big)$ are pairwise distinct and
pass through $F(d)$, that is $\big(F(d),f(A),f(B),f(C)\big)\in\Kk$.
Then $T$ is of full measure, and hence the set
\begin{equation*}
E'=
\big\{
e\in E:\,m_i(e)\in T,\,1\leq i\leq7
\big\}
\end{equation*}
is of full measure.
But it follows then from 
{\it la m\'ethode du quadrilat\`ere complet} \cite{Cartan}
that if $(d',d,A,B,C,M)\in E'$, and 
\begin{equation*}
\begin{aligned}
D\in d\hbox{ with }[A,B,C,D]=-1\\
D'\in F(d)\hbox{ with }\big[f(A),f(B),f(C),D'\big]=-1
\end{aligned}
\end{equation*}
then $f(D)=D'$.
\end{proof}

\begin{proof}[Continuation of the proof of Proposition~\ref{prop:cartan}]
Applying this argument to every $z\in E$,
we conclude that our original map $g$ has the property that
\begin{equation}\label{eq:g}
\begin{aligned}
\hbox{ for a.\,e. }z_1,z_2,z_3,z_4\in\CC\hbox{ such that }
[z_1,z_2,z_3,z_4]=-1,\\
\hbox{ we have that }
\big[g(z_1),g(z_2),g(z_3),g(z_4)\big]=-1\,.
\end{aligned}
\end{equation}
Fix now any $z_4$ such that (\ref{eq:g}) holds and 
consider the composition $\tilde g:=i_2\circ g\circ i_1$,
where $i_1,i_2$ are inversions with $i_1(\infty)=z_4$ and
$i_2(g(z_4))=\infty$.  
Then, for almost every $(z_1,z_2)\in\CC^2$, we have that
\begin{equation*}
2\tilde g\left(\frac{z_1+z_2}{2}\right)=\tilde g(z_1)+\tilde g(z_2)\,,
\end{equation*}
which implies that $\tilde g$ coincides almost everywhere
with an $\RR$-affine transformation of $\CC$. 
But since $\tilde g$ sends circles to circles, 
it is either $\CC$-affine or $\overline\CC$-affine,
so that $g$ is either a homography or an antihomography.
But then, using that $g$ has to preserve cyclic order 
on circles, one concludes that $g$ is $\CC$-affine.
\end{proof}

\begin{proof}[Proof of Theorem~\ref{thm:cartan}]
Let $\xi\in\partial\Hh_\CC^2$, $P$ the stabilizer of $\xi$ in 
$\mathrm{SU}(2,1)$ and $N$ its unipotent radical. 
Then $\partial\Hh_\CC^2\setminus\{\xi\}$ is a principal homogeneous
space for $N$ and the orbits of the center $Z(N)$ of $N$
correspond to the chains through $\xi$.  Fixing an identification
of $N/Z(N)$ with the additive group $\CC$ leads to a quotient map
\begin{equation*}
\pi_\xi:\partial\Hh_\CC^2\setminus\{\xi\}\twoheadrightarrow\CC
\end{equation*}
whose fibers are the chains through $\xi$ (with $\xi$ removed).
One can then fix an identification
\begin{equation*}
\xymatrix@1{P/Z(N)\ar[r]^\cong&\operatorname{Aff}(\CC)}
\end{equation*}
such that for the corresponding quotient homomorphism
\begin{equation*}
\omega_\xi:P\twoheadrightarrow\operatorname{Aff}(\CC)\,,
\end{equation*}
$\pi_\xi$ is equivariant with respect to $\omega_\xi$.
We recall here from \cite{Goldman_book} the following two facts:
\begin{enumerate}
\item[$(*)$] For every chain $C\subset\partial\Hh_\CC^2\setminus\{\xi\}$,
$\pi_\xi|_C$ is injective with image a circle in $\CC$;
\item[$(**)$] for every circle $S\subset\CC$ and any 
$s\in\partial\Hh_\CC^2\setminus\{\xi\}$ with $\pi_\xi(s)\in S$,
there is a (unique) chain $C$ through $s$ such that $\pi_\xi(c)=S$.
\end{enumerate}

Let us denote by $\Cc_p(\xi)$ the set of chains through the point
$\xi\in\partial\Hh_\CC^p$
and let now $\varphi:\partial\Hh_\CC^2\to\partial\Hh_\CC^2$
be a measurable map satisfying the hypotheses of Theorem~\ref{thm:cartan}.
Let $\Phi:\Cc_2\to\Cc_2$ be the map induced almost everywhere
on the set of chains defined in (\ref{eq:big-fhi}), 
and $E\subset\partial\Hh_\CC^2$ the subset of
full measure such that for every $\xi\in E$ and almost every $C\in\Cc_2(\xi)$,
also $\Phi(C)\in\Cc_2(\varphi(\xi))$.  Fix $\xi\in E$; composing
with an element from $\mathrm{SU}(2,1)$, we may assume that $\varphi(\xi)=\xi$.
Then $\Phi:\Cc_2(\xi)\to\Cc_2(\xi)$ induces a measurable map
$g_\xi:\CC\to\CC$ such that the diagram
\begin{equation*}
\xymatrix{
 \partial\Hh_\CC^2\setminus\{\xi\}\ar[r]^\varphi\ar[d]_{\pi_\xi}
&\partial\Hh_\CC^2\setminus\{\xi\}\ar[d]^{\pi_\xi}\\
 \CC\ar[r]^{g_\xi}
&\CC
}
\end{equation*}
commutes and $g_\xi$ satisfies the assumptions 
of Proposition~\ref{prop:cartan}, hence coincides almost everywhere
with an element of $\operatorname{Aff}(\CC)$.  
Thus composing $\varphi$ with an element $h\in P$ such that 
$\omega_\xi(h)=g_\xi$ almost everywhere, we may assume that 
$g_\xi=\id$ almost everywhere.  That is, for almost every $C\in\Cc_2(\xi)$ and almost
every $\xi\in C$, $\varphi(\xi)\in C$.  Now pick such a $C$ and $\eta\in C\cap E$.
Composing with an element from $Z(N)$, we may assume that $\varphi(\eta)=\eta$.
But then the map $g_\eta:\CC\to\CC$ fixes $\pi_\eta(\xi)$, leaves all circles 
through $\pi_\eta(\xi)$ invariant and coincides with an affine map
almost everywhere.  Hence $g_\eta=\id_\CC$ almost everywhere.
Then for almost every chain $C\in\Cc_2(\xi)$, $\pi_\eta$ being injective
on $C$ and $g_\eta=\id_\CC$, $\varphi|_C=\id_C$ and since those chains
foliate a set of $\partial\Hh_\CC^2$ of full measure, we conclude that 
$\varphi$ coincides almost everywhere with the identity.
\end{proof}

\vskip1cm
\section{Proof of Theorems~\ref{thm:thm6},~\ref{thm:thm7b} and~\ref{thm:thm8b}}\label{sec:proof}
\begin{proof}[Proof of Theorem~\ref{thm:thm6}]
Assume that $i_\rho=\mathrm{t}_\mathrm{b}(\rho)=1$
so that $\rho(\Gamma)\subset\mathrm{PU}(q,1)$ is nonelementary
and let $\varphi:\partial\Hh_\CC^p\to\partial\Hh_\CC^q$
be a $\Gamma$-equivariant measurable map.
Corollary~\ref{cor:formula} implies that
$\varphi$ satisfies condition (i) of Theorem~\ref{thm:cartan}.
Under the hypothesis that $p\geq2$, we show now that it satisfies
also condition (ii).  Assume that (ii) fails.
As before, let $\Cc_p(\xi)$ denote the set of chains through a point 
$\xi\in\partial\Hh_\CC^p$ and $\Phi:\Cc_p\to\Cc_q$ 
the measurable map defined in (\ref{eq:big-fhi}).  
Then for a set of positive measure of $\xi$'s, the set
\begin{equation*}
\big\{(C_1,C_2)\in\Cc_p(\xi)\times\Cc_p(\xi):\,\Phi(C_1)=\Phi(C_2)\big\}
\end{equation*}
is of positive measure.  In particular, for some $\xi\in\partial\Hh_\CC^p$
and $C_1\in\Cc_p(\xi)$, the set of $C_2\in\Cc_p(\xi)$
with $\Phi(C_1)=\Phi(C_2)$ is of positive measure.  
Now if $E\subset\Cc_p(\xi)$ is a set of positive measure, 
then the convex hull of $\cup_{C\in E}C$ is of full measure 
in $\partial\Hh_\CC^p$, which implies that 
the essential image $\essim(\varphi)$ of $\varphi$ 
is contained in the chain $\Phi(C_1)$ and hence that
$\rho(\Gamma)$ stabilizes a complex geodesic.  
We may thus assume that $q=1$.  

\begin{claim*} Let $\Ll_\rho\subset\partial\Hh_\CC^1=S^1$ be 
the limit set of $\rho(\Gamma)$.
Then $\Ll_\rho=S^1$.
\end{claim*}

\begin{proof} Observe first that since 
$\rho(\Gamma)\subset\mathrm{PU}(1,1)$
is nonelementary, it is Zariski dense, hence either discrete
or dense in $\mathrm{PU}(1,1)$.  If now $\Ll_\rho\neq S^1$,
then $\rho(\Gamma)$ is discrete and finitely generated,
and hence contains a normal torsionfree subgroup 
$\Lambda$ of finite index.  Since $\Ll_\rho=\Ll_\Lambda\neq S^1$,
the quotient $\Lambda\backslash\Hh_\CC^1$ is a noncompact
surface, hence $\Lambda$ is a free group, which implies
$\h^2(\rho(\Gamma),\RR)=\h^2(\Lambda,\RR)=0$,
and hence $\rho^{(2)}(\kappa_1)=0$ in $\h^2(\Gamma,\RR)$.
Since however $i_\rho=1$, clearly $\rho^{(2)}(\kappa_1)\neq0$,
which is a contradiction.  This shows that $\Ll_\rho=S^1$.
\end{proof}

\begin{claim*} With the above hypotheses $\ker\rho$ is finite and $\Gamma$ is 
cocompact.
\end{claim*}

\begin{proof}
Let us start by observing that if $I\subset S^1$ is any interval with nonvoid 
interior, then $\varphi^{-1}(I)$ contains, up to a null set,
an open subset of $\partial\Hh_\CC^p$.  

Indeed,
since $\essim(\varphi)=Sn^1$, then $\varphi^{-1}(I)$ is of positive
measure; moreover, for any $a\neq b$ in $\varphi^{-1}(I)$
such that $[\varphi(a),\varphi(b)]\subset I$,
the interval 
\begin{equation*}
\big\{z\in\partial\Hh_\CC^p:\,c_p(a,z,b)=1\big\}
\end{equation*}
belongs to $\varphi^{-1}(I)$ (up to measure zero),
which implies easily the assertion. 

To prove the claim, assume that $N:=\ker\rho<\Gamma$ is infinite.
Being discrete, its limit set in $\partial\Hh_\CC^p$ is nonvoid,
hence equals $\partial\Hh_\CC^p$, which implies that
$N$ acts minimally on $\partial\Hh_\CC^p$.  
Pick any interval $I\subset S^1$ with nonvoid interior
and let $\Oo\subset\partial\Hh_\CC^p$ be an open set such that 
$\Oo$ is included in $\varphi^{-1}(I)$ up to a set of measure zero.
Then $\varphi$ being $N$-invariant, and $N$ acting minimally
on $\partial\Hh_\CC^p$, we have that $\partial\Hh_\CC^p=\cup_{n\in N}n\Oo$
is contained in $\varphi^{-1}(I)$, up to measure zero.
But since $I$ was arbitrary, this is a contradiction.
This shows that $\ker\rho$ is finite.  If $\Gamma$ were not
cocompact, then -- since $p\geq2$ -- it would contain 
an integer Heisenberg group which would be sent, almost
injectively, into $\mathrm{PU}(1,1)$.  Since this is impossible,
it follows that $\Gamma$ is cocompact.  
\end{proof}

Thus $\Gamma$ and $\rho(\Gamma)$ are commensurable,
and hence their virtual cohomological dimensions coincide;
thus $\rho(\Gamma)$ has virtual cohomological dimension $4$
and hence cannot be discrete in $\mathrm{PU}(1,1)$.
Being Zariski dense, $\rho(\Gamma)$ is therefore dense in $\mathrm{PU}(1,1)$.
Passing to a subgroup of finite index of $\Gamma$, 
we may in addition assume that $\Gamma$ is torsionfree
and $\rho$ is injective.  Since the set of elliptic elements
is open in $\mathrm{PU}(1,1)$, we may pick $\gamma\neq \id$,
with $\rho(\gamma)$ elliptic.  Since $\Gamma$ is cocompact and torsionfree,
then $\gamma$ is necessarily hyperbolic.  Now pick a pair of open
intervals $\emptyset\neq I\subset I'$ such that 
the complement of $I'$ is of nonvoid interior.
Let $\Oo\subset\partial\Hh_\CC^p$ be nonvoid open subset 
such that $\Oo\subset\varphi^{-1}(I)$ up to a set of measure zero.
Conjugating by an element of $\Gamma$, we may assume
that the repulsive fixed point of $\gamma$ is in $\Oo$.
Let now $\{n_k\}_{k\in\NN}$ be a divergent sequence of integers
such that $\lim_{k\to\infty}\rho(\gamma)^{n_k}=\id$ in $\mathrm{PU}(1,1)$;
we may assume that $\rho(\gamma)^{n_k}I\subset I'$ for all $k\geq1$.
Then $\cup_{k\geq1}\gamma^{n_k}\Oo=\partial\Hh_\CC^p\setminus\{\xi\}$,
where $\xi$ is the attractive fixed point, and hence 
$\varphi^{-1}(I')$ equals $\partial\Hh_\CC^p$ up to a set
of measure zero.  Since $I'$ was arbitrary,
this is a contradiction.  

Thus we finally conclude that $\varphi$ satisfies also 
Theorem~\ref{thm:cartan}(ii), and hence there is a unique 
embedding $F:\Hh_\CC^p\to\Hh_\CC^q$ which is isometric and holomorphic
such that $\partial F=\varphi$ almost everywhere.
The uniqueness then implies that $F$ is $\Gamma$-equivariant.
\end{proof}

\begin{proof}[Proof of Theorem~\ref{thm:thm7b}]
By Corollary~\ref{cor:formula}, we have that
\begin{equation}\label{eq:cq=c1}
c_q\big(\varphi(\xi),\varphi(\eta),\varphi(\zeta)\big)=c_1(\xi,\eta,\zeta)
\end{equation}
for almost every $(\xi,\eta,\zeta)\in(\partial\Hh_\CC^1)^3$.  
Fix $\xi\neq\eta$ such that (\ref{eq:cq=c1}) holds for almost 
every $\zeta\in\partial\Hh_\CC^1$.  Then the essential image of $\varphi$ 
is contained in the chain $C$ determined by $\varphi(\xi)$ and
$\varphi(\eta)$, from which readily follows that $\rho(\Gamma)$
leaves invariant the complex geodesic whose boundary is $C$.
\end{proof}

\begin{proof}[Proof of Theorem~\ref{thm:thm8b}]
Let now $\rho:\Gamma\to\mathrm{PU}(1,1)$ be a homomorphism
with $|i_\rho|=1$ and, identifying $\partial\Hh_\CC^1=\SS^1$,
let  $\varphi:\SS^1\to\SS^1$ be the $\Gamma$-equivariant measurable map.
Then (\ref{eq:cq=c1}) holds and $\varphi$ is weakly order preserving,
so that \cite[Proposition~5.5]{Iozzi_ern} implies that there exists
a degree one monotone surjective continuous map $F:\SS^1\to\SS^1$
such that $f(\rho(\gamma)x)=\gamma f(x)$ for all $\gamma\in\Gamma$
and all $x\in\SS^1$.  The surjectivity of $f$ then implies that 
$\rho$ is injective, while its continuity that $\rho(\Gamma)$ is
discrete.  

Let now $\Ll\subset\SS^1$ be the limit set of $\rho(\Gamma)$.
To complete the proof we shall see that if the image under $\rho$ of a parabolic
element is also parabolic, then $\Ll=\SS^1$ and hence $\rho(\Gamma)$
is a lattice.  The proof will then be concluded by observing that 
if $\rho(\Gamma)$ is a lattice,
then it acts minimally on $\SS^1$ and hence $f$ is also injective.

Let us suppose by contradiction that $\Ll\subsetneq\SS^1$.
Then, since $\rho(\Gamma)\backslash\Hh_\CC^1$ is a complete hyperbolic 
surface topologically of finite type, 
for any connected component $I$ of $\SS^1\setminus\Ll$, one has that 
\begin{equation*}
\operatorname{Stab}_{\rho(\Gamma)}(I)=\,<\rho(\gamma)>\,,
\end{equation*}
where $\rho(\gamma)\in\Gamma$ is the hyperbolic element whose fixed points
are the endpoints $a,b$ of the interval $\overline I$.
Since $f$ is a semiconjugacy in the sense of Ghys \cite{Ghys_87},
the set $f(\overline I)$ is reduced to a point, say $f(\overline I)=\{\xi\}$.
Clearly $\gamma\xi=\xi$ and we shall show that this is the only 
fixed point of $\gamma$ hence $\gamma$ is parabolic contradicting the hypothesis.
Let us suppose that there exists $\eta\in\SS^1$ with $\gamma\eta=\eta$
and, since $f$ is surjective, let $x\in\SS^1$ such that $f(x)=\eta$.
Without loss of generality let us assume that $x\neq b$ 
(otherwise $\xi=\eta$ and we are done) and that $\rho(\gamma)^nx\to a$.
Then
\begin{equation*}
\eta=\gamma^n\eta=f\big(\rho(\gamma)^nx\big)\to f(a)=\xi\,,
\end{equation*}
and hence $\xi=\eta$ showing that $\gamma$ is parabolic.
\end{proof}

%%% Local Variables: 
%%% mode: latex
%%% TeX-master: "deform"
%%% End: 

\vskip1cm
\appendix
\section{Proof of Proposition~\ref{prop:resol}}

For the proof of Proposition~\ref{prop:resol}(1)
we need to show the existence of norm one contracting homotopy operators
from $\linfty\big((G/Q)^{n+1}_f\big)$ to $\linfty\big(\fib\big)$
sending $L$-continuous vectors into $L$-continuous vectors.

To this purpose we use the map $q_n$ which identifies 
the complex of Banach $G$-modules $\big(\linfty(G/Q)^\bullet_f\big)$ 
with the subcomplex $\big(\linfty(G\times(H/Q)^\bullet)^H\big)$
of $H$-invariant vectors of the complex $\big(\linfty(G\times(H/Q)^\bullet)\big)$,
where now the differential $d_n$ is given by 
\begin{equation*}
d_nf(g,x_1,\dots,x_n)=\sum_{i=0}^n(-1)^if(g,x_1,\dots,\hat{
x_i},\dots,x_n)\,,
\end{equation*}
and we show more generally that:

\begin{lemma}\label{lem:contr-hom} For every $n\geq0$ there are linear maps
\begin{equation*}
h_n:\linfty\big(G\times(H/Q)^{n+1}\big)\to\linfty\big(G\times(H/Q)^n\big)
\end{equation*}
such that:
\begin{enumerate}
\item $h_n$ is norm-decreasing and $H$-equivariant;
\item for any closed subgroup $L<G$, the map 
$h_n$ sends $L$-continuous vectors into $L$-continuous vectors, and
\item for every $n\geq1$ we have the identity
\begin{equation*}
h_nd_n+d_{n-1}h_{n-1}=\id\,.
\end{equation*}
\end{enumerate}
\end{lemma}

The Lemma~\ref{lem:contr-hom} and the remarks preceding it 
imply then Proposition~\ref{prop:resol}.

The construction of the homotopy operator in Lemma~\ref{lem:contr-hom}
requires the following two lemmas, 
the first of which showing that the measure $\nu$ on $H/Q$
can be chosen to satisfy certain regularity properties, 
and the second constructing an appropriate Bruhat function for $H<G$.

Let $dh$ and $d\xi$ be the left invariant Haar measures on $H$ and $Q$.

\begin{lemma}\label{lem:nu} There is an everywhere positive continuous function
$q:H\to\RR^+$ and  Borel probability measure $\nu$ on $H/Q$
such that
\begin{equation*}
\int_{H/Q}d\nu(x)\int_Qf(x\xi)d\xi=\int_Hf(h)q(h)dh\,,
\end{equation*}
for every $f\in\mathrm{C}_{00}(H)$.
\end{lemma}

\begin{proof} Let $q_1:H\to\RR^+$ be an everywhere positive
continuous function satisfying
\begin{equation*}
q_1(x\eta)=q_1(x)\frac{\Delta_Q(\eta)}{\Delta_H(\eta)}\,,\forall \eta\in Q\,\,x\in H\,,
\end{equation*}
where $\Delta_Q$, $\Delta_H$ are the respective modular functions
(see \cite{Reiter}), and let $\nu_1$ be the corresponding
positive Radon measure on $H/Q$ such that the above formula holds.
Then choose $q_2:H/Q\to\RR^+$ continuous and everywhere positive,
such that $q_2d\nu_1$ is a probability measure.  Then the
lemma holds with $q=q_1q_2$ and $\nu=q_2\otimes\nu_1$.
\end{proof}

A direct computation shows that 
\begin{equation}\label{eq:lambda}
\int_{H/Q}f(y^{-1}x)d\nu(x)=\int_{H/Q}f(x)\lambda_y(x)d\nu(x)\,,
\end{equation}
where $\lambda_y(x)=q(yx)/q(x)$, for all $f\in\mathrm{C}_{00}(H/Q)$ and $h\in H$.
In particular, the class of $\nu$ is $H$-invariant since $\lambda_y$
is continuous and everywhere positive on $H/Q$.

\begin{lemma}\label{lem:beta} There exists a function $\beta:G\to\RR^+$
such that
\begin{enumerate}
\item for every compact set $K\subset G$, $\beta$ coincides on $KH$
with a continuous function with compact support;
\item $\int_H\beta(gh)dh=1$ for all $g\in G$, and
\item $\lim_{g_0\to e}\sup_{g\in G}\int_H\big|\beta(g_0gh)-\beta(gh)\big|dh=0$
\end{enumerate}
\end{lemma}

\begin{proof} Let $\beta_0$ be any function satisfying (1) and (2) 
(see \cite{Reiter})
and let $f\in\mathrm{C}_{00}(G)$ be any nonnegative function
normalized so that 
\begin{equation*}
\int_Gf(x)d_rx=1\,,
\end{equation*}
where $d_rx$ is a right invariant Haar measure on $G$.
Define
\begin{equation*}
\beta(g)=\int_Gf(gx^{-1})\beta_0(x)d_rx\,,\qquad g\in G\,.
\end{equation*}
It is easy to verify that also $\beta$ satisfies (1) and (2),
and, moreover, it satisfies (3) as well.  
In fact, we have that for all $g_0,g\in G$, $h\in H$
\begin{equation*}
\beta(g_0gh)-\beta(gh)=\int_G\big(f(g_0gx^{-1})-f(gx^{-1})\big)\beta_0(xh)d_rx\,,
\end{equation*}
which implies, taking into account that $\int_G\beta_0(xh)dh=1$
and the invariance of $d_rx$, that 
\begin{equation*}
\int_H\big|\beta(g_0gh)-\beta(gh)\big|dh
\leq
\int_G\big|f(g_0x^{-1})-f(x^{-1})\big|d_rx\,,
\end{equation*}
so that
\begin{equation*}
    \lim_{g_0\to e}\sup_{g\in G}\int_H\big|\beta(g_0gh)-\beta(gh)\big|dh
\leq\lim_{g_0\to e}\int_G\big|f(g_0x^{-1})-f(x^{-1})\big|d_rx=0\,.
\end{equation*}
\end{proof}

\begin{proof}[Proof of Lemma~\ref{lem:contr-hom}]  
Let $\nu$ be as in Lemma~\ref{lem:nu} and $\beta$ as in Lemma~\ref{lem:beta}.
define a function
\begin{equation*}
\psi:G\times H/Q\to\RR^+
\end{equation*}
by 
\begin{equation*}
\psi(g,x):=\int_H\beta(gh)\lambda_{h^{-1}}(x)dh\,,
\end{equation*}
where $\lambda_h(x)$ is as in (\ref{eq:lambda}).
The following properties are then direct verifications:
\begin{enumerate}
\item $\psi(gh^{-1},hx)\lambda_h(x)=\psi(g,x)$ for all $g\in G$, $h\in H$ and
$x\in H/Q$;
\item $\int_{H/Q}\psi(g,x)d\nu(x)=1$, for all $g\in G$;
\item $\psi\geq0$ and is continuous.
\end{enumerate} 
This being, define for $n\geq0$ and $f\in \linfty(G\times(H/Q)^{n+1})$:
\begin{equation*}
h_nf(g,x_1,\dots,x_n)=\int_{H/Q}\psi(g,x)f(g,x_1,\dots,x_n,x)d\nu(x)\,.
\end{equation*}
Then, $h_nf\in\linfty\big(G\times(H/Q)^n\big)$ and (2) implies that
$\|h_nf\|_\infty\leq\|f\|_\infty$.  The fact that $h_n$ 
is an $H$-equivariant  homotopy operator is a formal
consequence of (1) and (2).

Finally, let $L<G$ be a closed subgroup and 
$f\in\linfty\big(G\times(H/Q)^{n+1}\big)$ an $L$-continuous vector, 
that is
\begin{equation*}
\lim_{l\to e}\|\theta(l)f-f\|_\infty=0\,,
\end{equation*}
where
\begin{equation*}
\big(\theta(l)f\big)(g,x_1,\dots,x_{n+1})=f(lg,x_1,\dots,x_n)\,.
\end{equation*}
Then
\begin{equation*}
\begin{aligned}
 &h_nf(lg,x_1,\dots,x_n)-h_nf(g,x_1,\dots,x_n)\\
=&\int_{H/Q}\psi(lg,x)\big(f(lg,x_1,\dots,x_n,x)-f(g,x_1,\dots,x_n,x)\big)d\nu(x)\\
+&\int_{H/Q}\big(\psi(lg,x)-\psi(g,x)\big)f(g,x_1,\dots,x_n,x)d\nu(x)\,.
\end{aligned}
\end{equation*}
The first term is bounded by $\|\theta(l)f-f\|_\infty$
taking into account (2), while the second is bounded by 
$\|f\|_\infty\int_{H/Q}\big(\psi(lg,x)-\psi(g,x)\big)d\nu(x)$.
Now
\begin{equation*}
\psi(lg,x)-\psi(g,x)=\int_H\big(\beta(lgh)-\beta(gh)\big)\lambda_{h^{-1}}(x)dh\,,
\end{equation*}
which, taking into account that 
$\int_{H/Q}\lambda_{h^{-1}}(x)d\nu(x)=1$,
implies that
\begin{equation*}
    \int_{H/Q}\big|\psi(lg,x)-\psi(g,x)\big|d\nu(x)
\leq\int_{H/Q}\big|\beta(lgh)-\beta(gh)\big|dh\,.
\end{equation*}
Thus
\begin{equation*}
\begin{aligned}
 \|\theta(l)h_nf-h_nf\|_\infty\leq
&\|\theta(l)f-f\|_\infty\\+
&\|f\|_\infty\sup_{g\in G}\int_H\big|\beta(lgh)-\beta(gh)\big|dh
\end{aligned}
\end{equation*}
which, using Lemma~\ref{lem:beta}, implies that
\begin{equation*}
\lim_{l\to e}\|\theta(l)h_nf-h_nf\|_\infty=0
\end{equation*}
and shows that $h_nf$ is an $L$-continuous vector.
\end{proof}

%%% Local Variables: 
%%% mode: latex
%%% TeX-master: "deform"
%%% End: 

\vskip1cm

\bibliographystyle{amsplain}
\bibliography{refs}
\vskip1cm

\providecommand{\bysame}{\leavevmode\hbox to3em{\hrulefill}\thinspace}
\providecommand{\MR}{\relax\ifhmode\unskip\space\fi MR }
\providecommand{\MRhref}[2]{%
  \href{http://www.ams.org/mathscinet-getitem?mr=#1}{#2}
}
\providecommand{\href}[2]{#2}

\vskip1cm
\end{document}